\documentclass{amsart}
\usepackage{amsmath,latexsym,amsfonts,amssymb,amsthm,mathrsfs,longtable,lscape,dsfont,fancybox,yfonts}
\usepackage{amsthm}
\usepackage{amssymb}
\usepackage{amscd}
\usepackage{amsfonts}
\usepackage{amsbsy}
\usepackage{fancyhdr}
\usepackage{graphicx}
\usepackage[all]{xy}
\usepackage{epsfig}
\usepackage{epstopdf}
\newtheorem {teo} {Theorem} [section]

\newtheorem {prop} [teo]{Proposition}
\newtheorem {cor} [teo] {Corollary}

\begin{document}
\title[Families of periodic orbits in the planar Hill's four-body problem]{Families of periodic orbits in the planar Hill's four-body problem}
\author[Jaime Burgos--Garc\'ia]{Jaime Burgos--Garc\'ia }
\address{Departamento Acad\'emico de Matem\'aticas \\
Instituto Tecnol\'ogico Aut\'onomo de M\'exico\\
01080 M\'exico, D. F.}%
\email{jbg84@xanum.uam.mx}%


\maketitle

\begin{abstract}
In this work we perform a numerical exploration of the families of planar periodic orbits in the Hill's approximation in the restricted four body problem, that is, after a symplectic scaling, two massive bodies are sent to infinity, by mean of expanding the potential as a power series in  $m_3^{1/3}$, (the mass of the third small primary) and taking the limit case when  $m_3\rightarrow 0$. The limiting Hamiltonian depends on a parameter $\mu$ (the mass of the second primary) and possesses some dynamical features from both the classical restricted three-body problem and the restricted four-body problem. We explore the families of periodic orbits of the infinitesimal particle for some values of the mass parameter, these explorations show interesting properties regarding the periodic orbits for this problem, in particular for the Sun-Jupiter-asteroid case. We also offer details on the horizontal and vertical stability of each family.
\end{abstract}

\section{Introduction}
For a long time the restricted three body problem (R3BP) which studies the dynamics of a massless particle moving under the gravitational force produced by two massive bodies in circular orbits around their center of mass, has been studied to know preliminary orbits in some space missions, it is worth mentioning works like \cite{Belbruno}, \cite{Broucke}, \cite{Koon} where the authors used this model to design trajectories of a spacecraft for some spatial missions. In recent years, the exploration of the Trojan asteroids was recognized by the 2013 Decadal Survey, which includes Trojan Tour and Rendezvous, among the New Frontiers missions in the decade 2013-2022. We refer by Trojan asteroid as an small body in a equilateral triangle configuration with other two bodies if they are viewed in a rotating frame of reference.  The most famous examples of Trojan asteroids in our solar system are  the Trojan asteroids of the Sun-Jupiter system, these groups of asteroids can be divided into two large groups, commonly named the `Trojans' and the `Greeks'. The first group is centered at a point on Jupiter's orbit around the Sun at $60^\circ$ behind the planet, and the second group is centered at a point  on the same  orbit at $60^\circ$ ahead the planet; therefore each of the two points forms an equilateral triangle configuration with the Sun and Jupiter. It is worth mentioning that there exist other Trojan asteroids in our solar system, for instance in the Mars-Sun and Neptune-Sun systems, Saturn and some of its moons form equilateral triangle configurations, Saturn--Tethys--Telesto, Saturn--Tethys--Calypso or Saturn--Dione--Helen are good examples of such configuration, therefore such triangle configurations are common in our solar system. Now, if we consider a small particle interacting with the three particles in this configuration, a four body approach becomes necessary. In several works like \cite{Moulton} \cite{PapaIII} ,\cite{Cecc}, \cite{Burgos} and \cite{BurgosII} the authors have studied the dynamics of the infinitesimal particle, the so called restricted four body problem (R4BP). However, in this paper we consider the Hill's approximation approach in the restricted four-body problem, such approximation was introduced by J. Burgos-Garcia and M. Gidea in \cite{BurgosIII}. This new approach describes the dynamics of a massless particle in a neighborhood of a small mass at one of the vertices of a Lagrange central configuration. The system produced, hereafter refereed as Hill four body problem (H4BP), provides an approximation of the motion of the massless particle in an $O(m_3^{1/3})$-neighborhood  of $m_3$, the remaining masses $m_1$ and $m_2$ are sent at infinite distance through the rescaling, furthermore, if we let  $m_2\to 0$ in our model, the resulting system coincides with the classical lunar Hill problem, so the H4BP is a generalization in this sense. In \cite{BurgosIII} the authors pointed out the main differences between the classical and the generalized problem, in particular the existence of two additional equilibrium points, the so called $L_{3}$ and $L_{4}$. These points have the property that there exists a value of the mass parameter $\mu_{0}$ for which such equilibrium points are linearly stable if $\mu<\mu_{0}$ and they are unstable for $\mu>\mu_{0}$, when $\mu=\mu_{0}$ we have two pure imaginary eigenvalues of multiplicity two. This property will have remarkable consequences in the exploration of the periodic orbits which is the main objective of this work.

Between the possible applications of the planar periodic orbits of this system, we can consider the families of these periodic orbits as starting point for the design of spacecraft trajectories near a Trojan asteroid of the solar system, or as preliminary orbits for modelling the dynamics of the particles in the so called E-ring of Saturn  where the systems Saturn--Tethys--Telesto, Saturn--Tethys--Calypso are contained \cite{Esposito}. Therefore, in order to obtain a first insight of the dynamics of the massless particle in such systems, an exploration of the families of periodic orbits of the system is required.

This paper is organized as follows: In section 2 we introduce the Hill four body problem, we state the equations of motion, their properties and the Hamiltonian of the system. In section 3 we recall the results regarding the equilibrium points of the system, however, we offer a deep study of the linear dynamics around the new equilibrium points produced by the presence of the second massive body, as we will see, such equilibrium points are linearly stable and they provide new families of periodic orbits that are not present in the classical Hill's problem. In section 4 we discuss the so called Levi-Civita regularization of the system, such regularization has been already discussed in \cite{BurgosIII}, however, we offer more details regarding the regularized equations, this study of the regularized equations will prove to be useful in the numerical computations. In section 5 we show the numerical explorations of the families of periodic orbits in the H4BP, we show the geometry of the orbits, their qualitative behaviour respect to the variations of the Jacobi constant and the evolution of their stability.

\section{The Hill´s approximation in the four body problem}

In \cite{BurgosIII} the authors introduced the H4BP by studying the limit when $m_{3}\rightarrow0$ in the Hamiltonian of the R4BP. The procedure was similar to that in \cite{MeyerS,MeyerHDS}, first, a symplectic scaling depending on $m_3^{1/3}$ was performed in the Hamiltonian, second, an expansion in the gravitational potential of the Hamiltonian as a power series in $m_3 ^{1/3}$ in a neighborhood of the small mass $m_{3}$, then studying the behavior of the Hamiltonian in the limit as $m_3\to 0$. The result of the cited procedure produces the following theorem whose proof can be found in the page 6 of \cite{BurgosIII}.

\begin{teo} \label{main theorem} After the symplectic scaling
$$(x,y,z,p_{x},p_{y},p_{z})\rightarrow m_{3}^{1/3}(x,y,z,p_{x},p_{y},p_{z}),
$$ the limit $m_{3}\rightarrow0$ of the Hamiltonian (\ref{originalhamiltonian}) restricted to a neighborhood of $m_{3}$ exists and yields a new Hamiltonian
\begin{equation}\begin{split}\label{hillhamiltonian}
H=&\frac{1}{2}(p^{2}_{x}+p^{2}_{y}+p_{z}^{2})+yp_{x}-xp_{y}+\frac{1}{8}x^2-\frac{3\sqrt{3}}{4}(1-2\mu)xy-
\frac{5}{8}y^2+\frac{1}{2}z^2\\
&-\frac{1}{\sqrt{x^2+y^2+z^2}},
\end{split}\end{equation}where $m_{1}=1-\mu$ and $m_{2}=\mu$.
\end{teo}

We note that the resulting Hamiltonian is a three-degree of freedom system depending on a parameter $\mu$ which becomes equal to the mass of the secondary~$m_{2}$. As it was pointed out in the introduction, this systems represents a Kepler problem, for the infinitesimal mass and the so called asteroid, plus quadratic terms produced by the gravitational influence of the two large bodies. In this sense, the current system represents a gravitational generalization of the classical Hill's problem, i.e., when we consider $\mu=0$ in Hamiltonian (\ref{hillhamiltonian}), we recover the system of the Hill's problem. This is not obvious at first sight, however if we apply a rotation in the $xy$-plane we obtain a new Hamiltonian with nicer properties, in particular, in these coordinates we recover the original Hill's problem. 

\begin{cor} \label{main corollary} The equations of motion given by (\ref{hillhamiltonian}) are equivalent, via a rotation, with the system
\begin{equation}\begin{split}\label{finalequations}
\ddot{\bar{x}}-2\dot{\bar{y}}&=\Omega_{\bar{x}},\\
\ddot{\bar{y}}+2\dot{\bar{x}}&=\Omega_{\bar{y}},\\
\ddot{\bar{z}}&=\Omega_{\bar{z}},\end{split}\end{equation}
with
\begin{equation}\label{rotatedeffective}
\bar{\Omega}=\frac{1}{2}(\lambda_{2}\bar{x}^2+\lambda_{1}\bar{y}^2-\bar{z}^2)+
\frac{1}{\sqrt{\bar{x}^2+\bar{y}^2+\bar{z}^2}},
\end{equation}
where $\lambda_{2}=\frac{3}{2}(1+d)$ and $\lambda_{1}=\frac{3}{2}(1-d)$ are the eigenvalues corresponding to the rotation transformation in the $xy$-plane with $d=\sqrt{1-3\mu+3\mu^2}$.
\end{cor}

The system possesses a Jacobi type first integral given by the expression
\begin{equation}\label{jacobiintegral}
\dot{x}^2+\dot{y}^2=2\Omega-C,
\end{equation}

where $C$ is a constant. The expressions for $\Omega_{x}$ and $\Omega_{y}$ have the following properties
\begin{eqnarray*}
\Omega_{x}(x,-y)&=&  \Omega_{x}(x,y),\\
\Omega_{y}(x,-y)&=&- \Omega_{x}(x,y).
\end{eqnarray*}
These properties state that the equations (\ref{finalequations}) are invariant under the transformations $x \rightarrow x$, $y \rightarrow -y$, $\dot{x} \rightarrow -\dot{x}$, $\dot{y} \rightarrow \dot{y}$, $\ddot{x} \rightarrow \ddot{x}$, $\ddot{y} \rightarrow -\ddot{y}$ as a consequence we have the well known symmetry respect the $x$-axis. A similar argument shows that the equations (\ref{finalequations}) are also symmetric respect the $y$-axis. Therefore, we have the well known symmetries of the classical Hill's problem that can be recovered by considering $\mu=0$, in such a case, the equations (\ref{finalequations}) take the form
\begin{equation*}\begin{split}
\ddot{\bar{x}}-2\dot{\bar{y}}&=\bar{\Omega}_{\bar{x}},\\
\ddot{\bar{y}}+2\dot{\bar{x}}&=\bar{\Omega}_{\bar{y}},\end{split}\end{equation*}
with $$\bar{\Omega}=\frac{3}{2}\bar{x}^2+\frac{1}{\Vert w\Vert}$$ exactly as in the classical Hill problem, see \cite{Sz}.

The Hamiltonian in these new coordinates is given by the expression
\begin{equation}\begin{split}
H(x,y,z,p_x,p_y,p_z)=&\frac{1}{2}(p_x^2+p_y^2+p_z^2)+yp_x-xp_y+ax^2+by^2+cz^2\\&-\frac{1}{\sqrt{x^2+y^2+z^2}},
\label{rotatedspatialhamiltonian}
\end{split}
\end{equation}
where $a=(1-\lambda_{2})/2$, $b=(1-\lambda_{1})/2$ and $c=1/2$.

Throughout this work we will restrict ourselves to the study of the periodic orbits for the planar case, $z=0$. In forthcoming work, we hope to provide more details on the three dimensional dynamics of this system by considering the information of the planar periodic orbits that is the main objective of this paper.

\section{The equilibrium points of the system.}

The previous work performed in \cite{BurgosIII} shows that the H4BP possesses four planar equilibrium points given by the expressions
\begin{eqnarray*}
L_{1}=\left(\frac{1}{\sqrt[3]{\lambda_{2}}},0\right),
L_{2}=\left(-\frac{1}{\sqrt[3]{\lambda_{2}}},0\right),
L_{3}=\left(0,\frac{1}{\sqrt[3]{\lambda_{1}}}\right),
L_{4}=\left(0,-\frac{1}{\sqrt[3]{\lambda_{1}}}\right),
\end{eqnarray*}

The linear stability of the equilibrium points $L_{1}$ and $L_{2}$ is provided by the following

\begin{prop} The coefficient $B$ is negative for $\mu\in[0,1/2]$ so the equilibrium points $L_{1}$ and $L_{2}$ are unstable for this range of values of the mass parameter, in fact, the eigenvalues are given by $\pm\Lambda$ and $\pm\textit{i}\omega$ with $\Lambda>0$ and  $\omega>0$.
\end{prop}

On the other hand, for the points $L_{3}$ and $L_{4}$ we have the following

\begin{prop} There exists a value $\mu_{0}$ such that $D=0$, as a consequence, the equilibrium points $L_{3}$ and $L_{4}$ have the following properties: for $\mu\in(0,\mu_{0})$ their eigenvalues are $\pm\textit{i}\omega_{1}$ and $\pm\textit{i}\omega_{2}$, for $\mu=\mu_{0}$ we have a pair of the eigenvalues $\pm\textit{i}\omega$ of multiplicity 2, finally when $\mu\in(\mu_{0},1/2]$ the eigenvalues are $\pm\alpha\pm\textit{i}\omega$ with $\alpha>0$ and $\omega>0$.
\end{prop}

The value where the change of stability occurs is
$$
\mu_{0}=\frac{1}{450}\left(225-\sqrt{3(5227+2368\sqrt{21})}\right),
$$
which is approximately $\mu_{0}\approx 0.011942$. In \cite{BurgosIII} is reported the value $\mu_{0}\approx 0.00898964$ which is slightly lesser than the corrected value reported in this work, however, it is worth mentioning that the conclusions of the above cited work remain true, in particular, for the solar system $\mu\in[0,0.00095]$ as a consequence the equilibrium points $L_{3}$ and $L_{4}$ are always linearly stable, therefore, for the four equilibrium points, the Liapunov center theorem assures that under non resonant conditions, families of periodic orbits emanate from these equilibrium points. For the equilibrium points $L_{3}$ and $L_{4}$we have two families of such orbits, the so called short and long period families. In order to consider these families of periodic orbits for a further analytic continuation, a more detailed discussion on the linear dynamics around these points is required. We recall that because of the symmetry respect the $y-$axis, it will enough to study one of these points, say, $L_{3}$. \\

\subsection{The linear dynamics around $L_{3}$}\label{lineardynamics}

In the following we will focus on the study the linear system 
\begin{equation}\label{linearmatrixsystem}
\boldmath{\dot{\bar{x}}}=\bf{A}\bf{\bar{x}},
\end{equation}

where $\bf{\bar{x}}\rm=(\xi,\eta,\dot{\xi},\dot{\eta})^{T}$ and $A$ is the matrix

\begin{equation}
\left(\begin{array}{cccc}0 & 0 & 1 & 0 \\0 & 0 & 0 & 1 \\\Omega_{xx} & \Omega_{xy} & 0 & 2 \\\Omega_{xy} & \Omega_{yy} & -2 & 0\end{array}\right).\label{matrixlinearization}
\end{equation}
The partial derivatives in the point $L_{3}$ are given by the expressions

 \begin{equation*}\begin{split}\Omega_{xx}&=\lambda_{2}-\lambda_{1},\\ \Omega_{yy}&=3\lambda_{1},\\ \Omega_{xy}&=0.\end{split}\end{equation*}

We recall that the spectrum of the points $L_{3}$ and $L_{4}$ is given by the expression 

$$\Lambda_{1,2,3,4.}=\pm\frac{1}{\sqrt{2}}\sqrt{-A\pm\sqrt{D}},$$

where $A=(3d-1)/2$ and $D=(225d^2-222d+1)/4$ with $d=\sqrt{3\mu^2-3\mu+1}$. For $\mu\in(0,\mu_{0}]$ the spectrum is ${\pm\textit{i}\omega_{1},\pm\textit{i}\omega_{2}}$ with

$$\omega_{1}=\frac{1}{\sqrt{2}}\sqrt{A-\sqrt{D}},$$
$$\omega_{2}=\frac{1}{\sqrt{2}}\sqrt{A+\sqrt{D}}.$$

It is easy to see that when $\mu=\mu_{0}$ we have $D=0$, therefore $\omega_{1}=\omega_{2}=\sqrt{\frac{A}{2}}$. A straightforward computation shows that 

\begin{equation}
\label{expressionsfrequences}
\omega_{1}^2+\omega_{1}^2 =A,
\end{equation}

$$0<\omega_{1}\leq\sqrt{\frac{A}{2}}\leq\omega_{2}$$
$$\omega_{1}\omega_{2} =\frac{1}{2}\sqrt{A^2-D},$$
$$\frac{\omega_{2}}{\omega_{1}} =\sqrt{\frac{A+\sqrt{D}}{A-\sqrt{D}}}.$$

\begin{figure}
\centering
\includegraphics[width=2.5in]{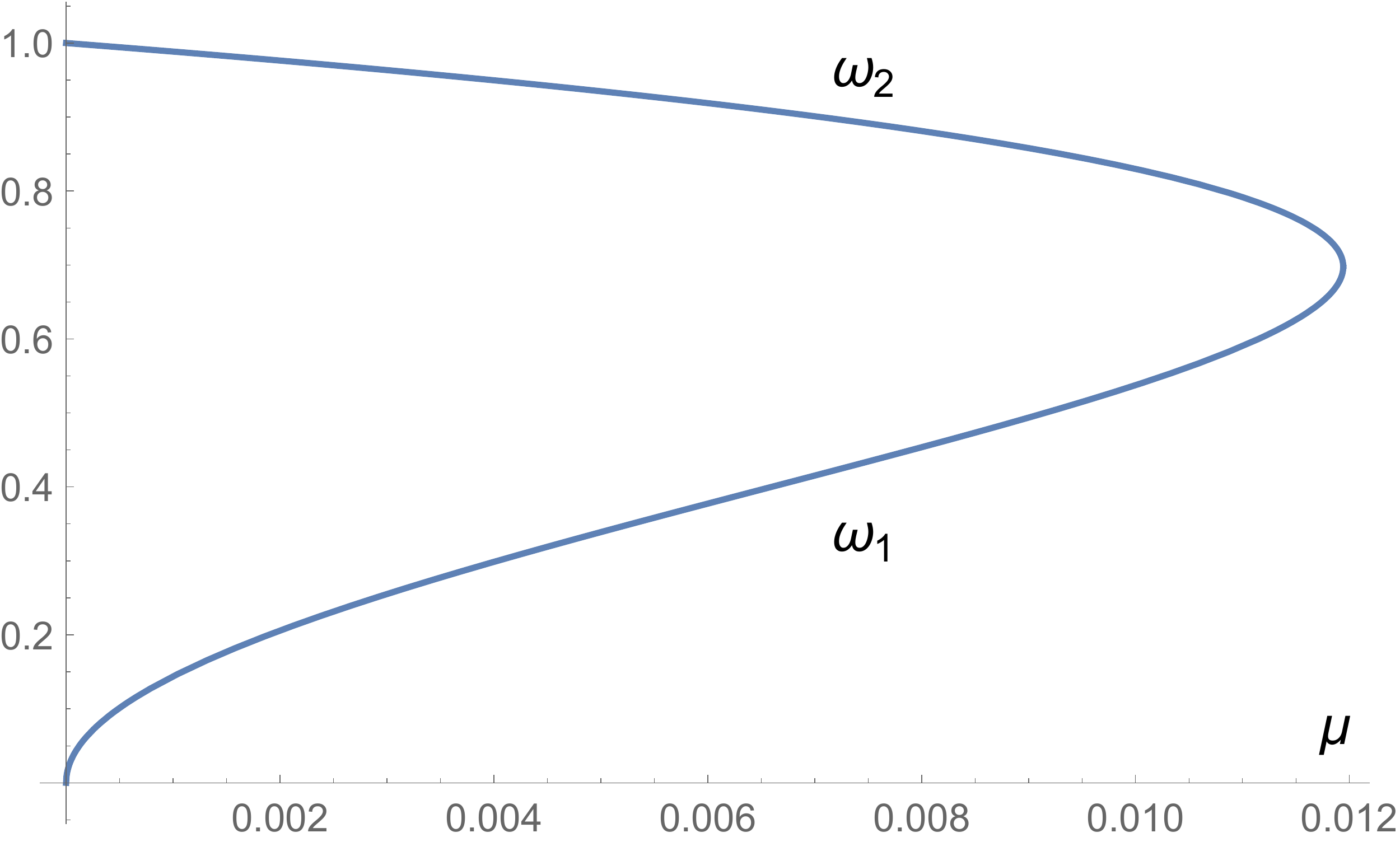}\label{frequencies}
\caption{Frequencies as functions of the mass parameter.\label{triangle}}
\end{figure}

The system (\ref{linearmatrixsystem}) can be written as

\begin{equation*}\begin{split}\ddot{\xi}-2\dot{\eta}&=\bar{\Omega}_{\xi\xi}\xi,\\ \ddot{\eta}+2\dot{\xi}&=\bar{\Omega}_{\eta\eta}\eta,\end{split}
\end{equation*}

with $\bar{\Omega}=\frac{1}{2}((\lambda_{2}-\lambda_{1})\xi^2+3\lambda_{1}\eta^2)$. It is clear that the $\xi-$axis is parallel to the original $x-$axis, therefore no rotation of coordinates is involved between the two coordinate systems.

Therefore, the general real solution of the linear system is given by
\begin{equation}\label{linearsolution}
\xi=K_{1}\cos\omega_{1}t+K_{2}\sin\omega_{1}t+K_{3}\cos\omega_{2}t+K_{4}\sin\omega_{2}t,
\end{equation}

$$\eta=\bar{K}_{1}\cos\omega_{1}t+\bar{K}_{2}\sin\omega_{1}t+\bar{K}_{3}\cos\omega_{2}t+\bar{K}_{4}\sin\omega_{2}t,$$

with $K_{i}$ and $\bar{K}_{i}$ constants that depend on the initial conditions for $i=1,2,3,4.$, these constants are not independent but they are related by 

\begin{equation}
\begin{split}
K_{1}&=-\frac{2\omega_{1}}{\lambda_{2}-\lambda_{1}+\omega_{1}^2}\bar{K}_{2},\\ K_{2}&=\frac{2\omega_{1}}{\lambda_{2}-\lambda_{1}+\omega_{1}^2}\bar{K}_{1},\\ K_{3}&=-\frac{2\omega_{2}}{\lambda_{2}-\lambda_{1}+\omega_{2}^2}\bar{K}_{4}, \\ K_{4}&=\frac{2\omega_{2}}{\lambda_{2}-\lambda_{1}+\omega_{2}^2}\bar{K}_{3}.\end{split}
\end{equation}

The relation between such constants and the initial conditions for $t=0$ is

\begin{equation}
\begin{split}
\xi_{0}&=K_{1}+K_{3},\\ \eta_{0}&=\bar{K}_{1}+\bar{K}_{3},\\ \dot{\xi}_{0}&=\omega_{1}K_{2}+\omega_{2}K_{4}, \\ \dot{\eta}_{0}&=\omega_{1}\bar{K}_{2}+\omega_{2}\bar{K}_{4}.\end{split}
\end{equation}

In the solution (\ref{linearsolution}) we have two kinds of motion, the so called short and long period solutions. The short period solution is determined by the choice of constants $K_{i}=\bar{K_{i}}=0$ for $i=1,2.$ and the long period solution is determined by $K_{i}=\bar{K_{i}}=0$ for $i=3,4.$ The initial conditions for the short period solution become $\xi_{0}=K_{3}$, $\eta_{0}=\bar{K}_{3}$, $\dot{\xi}_{0}=\omega_{2}K_{4}$ and $\dot{\eta}_{0}=\omega_{2}\bar{K}_{4}$. If we define $$\alpha_{1}=\frac{2\omega_{1}}{\lambda_{2}-\lambda_{1}+\omega_{1}^2},$$ and $$\alpha_{2}=\frac{2\omega_{2}}{\lambda_{2}-\lambda_{1}+\omega_{2}^2},$$ we can write the initial conditions for the short period as $\xi_{0}=K_{3}$, $\eta_{0}=\bar{K}_{3}$, $\dot{\xi}_{0}=\omega_{2}\alpha_{2}\eta_{0}$ and $\dot{\eta}_{0}=-\frac{\omega_{2}}{\alpha_{2}}\xi_{0}$. Therefore the short period solution looks like 

\begin{equation}
\begin{split}
\xi&=\xi_{0}\cos\omega_{2}t+\alpha_{2}\eta_{0}\sin\omega_{2}t,\\ \eta&=\eta_{0}\cos\omega_{2}t-\frac{\eta_{0}}{\alpha_{2}}\sin\omega_{2}t.
\end{split}
\end{equation}

This solution will be completely determined by the constants $\xi_{0}$ and $\eta_{0}$, in particular we note that the equation $\dot{\eta}_{0}=-\frac{\omega_{2}}{\alpha_{2}}\xi_{0}$ states that the short period motion is retrograde because $\omega_{2}$ and $\alpha_{2}$ are positive. The above equation can be rewritten in order to figure out its geometry, the new expression looks like 

\begin{equation}
\begin{split}
\xi&=A\cos(\omega_{2}t+\phi),\\ \eta&=B\cos(\omega_{2}t+\phi),
\end{split}
\end{equation}

with 
$$A=\sqrt{\xi_{0}^2+\alpha_{2}^2\eta_{0}^2},$$
$$B=\sqrt{\eta_{0}^2+\frac{\xi_{0}^2}{\alpha_{2}^2}},$$
and
$$\tan\phi=\frac{\alpha_{2}\eta_{0}}{\xi_{0}}.$$

The constants $A$ and $B$ satisfy $A=\vert\alpha_{2}\vert B$, a straightforward computation shows that $\alpha_{2}<1$, so $A<B$ and the major axis of the elliptic motion is on the $\eta-$axis. The eccentricity of the ellipse is given by

$$e=\frac{\sqrt{B^2-A^2}}{B}=\sqrt{1+\alpha_{2}^2},$$

so the eccentricity of the ellipses depend only on the mass parameter $\mu$. A similar study can be performed for the long period motion to obtain that this kind of motion is elliptic with major axis on the $\eta-$axis. For the particular value of the Sun-Jupiter system $\mu\approx 0.0095$ the value of the periods is $T_{s}\approx6.35271$ and $T_{l}=44.8422$ for the short and long period respectively.

In order to perform a numerical continuation from the infinitesimal orbits for each value of $\mu$ we need to consider the values for which there exists resonances between the short and long period motions, i.e., values of $\mu$ such that 

$$\frac{\omega_{2}}{\omega_{1}}=k,$$ 

where $k\in\mathbb{Z}$. The above equation implies that $\omega_{2}^2=k^2\omega_{1}^2$, and from the equation (\ref{expressionsfrequences}) we obtain $\omega_{1}^2=\frac{\alpha}{2}-\omega_{1},$ with $\alpha=3d-1$. Therefore we can write $\omega_{2}^2(1+k^2)=\frac{\alpha k^2}{2}$ or 
\begin{equation}\label{omegamu}
\omega_{2}^2=c\alpha,
\end{equation}

where $c=\frac{1}{2}\left(\frac{k^2}{k^2+1}\right)$. The discriminant $D$, which is involved in the expressions of the frequencies, can be written as $D=\sqrt{\beta}/2$ with $\beta=225d^2-222d+1$, therefore the equation (\ref{omegamu}) can be written as $\beta=\alpha^2(4c-1)^2$. Finally, the equations to be solved for $\mu$ as function of $k$ is 
\begin{equation}\label{quad1}
a(K)d^2+b(K)d+c(K)=0,
\end{equation}
with $a(K)=225-9K$, $b(K)=6K-222$, $c(K)=1-K$ and $K=(\frac{k^2-1}{k^2+1})^2$. It is easy to see that $0\leq K<1$ for $k=\pm1,\pm2,\pm3,....$ Therefore $a(K)>0$, $b(K)<0$ and $c(K)>0$. By solving two quadratic equations, one of them given by (\ref{quad1}) and the second one given by $d^2=3\mu^2-3\mu+1$ we obtain

\begin{equation}\label{resonancesequation}
\mu=\frac{1}{2}-\frac{1}{6\sqrt{3}}\sqrt{\frac{5227+1184\sqrt{84-3K}-5K^2-32K\sqrt{84-3K}-38K}{(K-25)^2}}.
\end{equation}

In the figures (\ref{frequencies}) we can observe the behaviour of the frequencies as a functions of the mass parameter, in the table (\ref{table:nonlin}) we provide some values of the parameter where we have resonances between the frequencies. It is worth mentioning that in the table we show the values of $\mu$ with 6 decimals, however the formula (\ref{resonancesequation}) can provide exact results, for instance the true values for $k=2,3$ are respectively

\begin{eqnarray*}
\mu=\frac{1}{2}-\frac{1}{462}\sqrt{\frac{5}{3}(10181+458\sqrt{2073})},
\mu=\frac{1}{2}-\frac{1}{1218}\sqrt{5(24077+6464\sqrt{57})}.
\end{eqnarray*}

\begin{table}[ht]
\caption{Resonances between the frequencies} 
\centering 
\begin{tabular}{c c c c}
\hline\hline 
k & $\mu$ & k & $\mu$ \\ [0.5ex] 
\hline 
1 & $\mu_{0}$ & 6 & 0.001293 \\ 
2 & 0.007733 & 7 & 0.000965 \\
3 & 0.004390 & 8 & 0.000746 \\
4 & 0.002713 & 9 & 0.000594 \\
5 & 0.001817 & 10 & 0.000483 \\ [1ex] 
\hline 
\end{tabular}
\label{table:nonlin} 
\end{table}

\section{The regularized equations}

In the numerical explorations we will find several orbits colliding with the so called tertiary, so in order to compute these kind of orbits and follow their evolution as the Jacobi constant is varied, we will need to perform the well known Levi-Civita regularization in the Hamiltonian of the system. A first study of the regularized system can be found again in \cite{BurgosIII}, however, in this work we will adopt a slightly different approach in the regularization process, this approach will show to be valuable for our computation purposes. The Hamiltonian of the planar system \eqref{rotatedspatialhamiltonian} is

\[H(x,y,p_x,p_y)=\frac{1}{2}\Vert p\Vert^2-q^{T}Kp+\frac{1}{2}q^{T}Dq-\frac{1}{\Vert q\Vert},\]

with $q=(x,y)$, $p=(p_{x},p_{y})$, $$K=\left(
\begin{array}{rr}  0 & 1 \\  -1 &  0 \\\end{array} \right),$$ and $$D=\left(
\begin{array}{rr}  1-\lambda_{2} & 0 \\  0 &  1-\lambda_{1} \\\end{array} \right).$$

The Levi-Civita procedure  consists in changing the coordinates and  the conjugate momenta and in rescaling the time, as follows:
\[\left(\begin{array}{c} x \\ y \\ \end{array} \right)\longrightarrow A_{0}
\left(\begin{array}{c}  Q_{1} \\  Q_{2} \\ \end{array} \right),\]
\[\left(\begin{array}{c} p_x \\ p_y \\ \end{array} \right)\longrightarrow\frac{1}{ 2(Q_{1}^2+ Q_{2}^2)}A_{0}
\left(\begin{array}{c}  P_{1} \\  P_{2} \\ \end{array} \right),\]
and \[d\tau\longrightarrow\frac{1}{ 4(Q_{1}^2+ Q_{2}^2)}dt,\]

with 
$$
A_{0}=\left(
\begin{array}{rr}  Q_{1} & - Q_{2} \\  Q_{2} &  Q_{1} \\\end{array} \right).
$$
The transformed Hamiltonian is \[\bar{H}(Q_{1},Q_{2},P_{1},P_{2})=4(Q_{1}^2+ Q_{2}^2)\left(H(Q_{1},Q_{2},P_{1},P_{2})+\frac{C}{2}\right),\] where $C$ is the Jacobi constant, and  $(Q_{1},Q_{2},P_{1},P_{2})$ denote the transformed variables.

A straightforward computation shows that the regularized Hamiltonian is given by the following expression:

\begin{equation}
\bar{H}=\frac{1}{2}\Vert P\Vert^2-2\Vert Q\Vert^2(P_{2}Q_{1}-Q_{2}P_{1})+2\Vert Q\Vert^2Q^{T}M_{2}Q+2C\Vert Q\Vert^2-4,
\label{eqn:hillreg1}
\end{equation}

where $Q=(Q_{1},Q_{2})$, $P=(P_{1},P_{2})$ and $M_{2}=A_{0}^{T}DA_{0}$. We can omit from $\hat H$ the constant $4$ because it contributes nothing to the equations of motion. If we eliminate the $P$ coordinate in the equations of the Hamiltonian (\ref{eqn:hillreg1}), they can be rewritten in terms of the $Q$ coordinates only, such equations looks like

\begin{equation}
\begin{split}\label{regularizedequations}
\ddot{Q_{1}}-8(Q_{1}^2+Q_{2}^2)\dot{Q_{2}^2}&=\Omega_{Q_{1}}^{r},\\ \ddot{Q_{2}}+8(Q_{1}^2+Q_{2}^2)\dot{Q_{1}^2}&=\Omega_{Q_{2}}^{r},
\end{split}
\end{equation}
 
where 
$$
\Omega_{Q_{1}}^{r}=-4Q_{1}(C-4(2Q_{1}^2Q_{2}^2+Q_{2}^4)\lambda_{1}+(2Q_{1}^2Q_{2}^2-3Q_{1}^4+Q_{2}^4)\lambda_{2}),
$$

$$
\Omega_{Q_{2}}^{r}=-4Q_{2}(C-4(2Q_{1}^2Q_{2}^2+Q_{1}^4)\lambda_{1}+(2Q_{1}^2Q_{2}^2-3Q_{2}^4+Q_{1}^4)\lambda_{2}),
$$

$$
\Omega^{r}=2(Q_{1}^2+Q_{2}^2)(4Q_{1}^2Q_{2}^2\lambda_{1}+(Q_{1}^2-Q_{2}^2)^2\lambda_{2})-2C(Q_{1}^2+Q_{2}^2)+4,
$$

the constant term $4$ can be omitted again in the above expression. The new effective regularized potential can be written as follows:

$$
\Omega^{r}=4(Q_{1}^2+Q_{2}^2)\left(\Omega(Q_{1},Q_{2})-\frac{C}{2}\right),
$$

but in this case $\Omega(Q_{1},Q_{2})$ is the physical effective potential evaluated in the regularized variables. The Jacobi first integral in the regularized variables looks like

$$
\dot{Q}_{1}^2+\dot{Q}_{2}^2=2\Omega^{r}.
$$

It is worth mentioning that the equations (\ref{regularizedequations}) are computationally cheaper than the ones of the Hamiltonian (\ref{eqn:hillreg1}), this fact will represent a great advantage in the time of computation.

\section{Numerical explorations}
In the present section we offer a numerical exploration of the families of periodic orbits for the planar H4BP, although the theory predicts that there are infinitely many families of periodic orbits, we restrict our explorations to the basic families of symmetric and asymmetric periodic orbits that intersect the $x-$axis two and four times, these orbits are commonly called \textit{simple-periodic} and \textit{double-periodic} orbits in the sense of E. Str\"{o}mgren \cite{Henon2003}. It is well known that the exploration of these basic families provide the backbone and a first insight on the families of periodic orbits with greater number of intersections, furthermore, if we consider the information provided by the so called vertical stability, we can obtain valuable information for the computation of spatial periodic orbits, for instance, halo orbits near the equilibrium points. See \cite{Gomez} for a good exposition on the dynamics around the collinear equilibrium point of the R3BP.

Due to the existence of the Jacobi type first integral (\ref{jacobiintegral}), the families of periodic orbits lie in a smooth cylinder parametrized by the Jacobi constant or equivalently, parametrized by the energy so, each family can be continued continuously respect the Jacobi constant $C$. Throughout this section we use the classical notation of M. H\'enon \cite{Henon}, \cite{Henon2003} to denote the families of periodic orbits. We use the classical $(x,C)$ plane of characteristic curves to represent the families of periodic orbits, in addition we use the $(a_{h},C)$ plane to show the evolution of the horizontal stability of the families, here $a_{h}=a+d$ denotes the horizontal stability index (see \cite{HenII} for details), we have stability in the linear sense when $\vert a_{h}\vert<2$. In the case of symmetric periodic orbits we have that $a=d$ and the stability criterion reduces to $\vert a_{h}\vert<1$. We also have computed the so called vertical stability index $a_{v}=a_{33}+a_{66}$, where $a_{33}$ and $a_{66}$ are the respective elements of the monodromy matrix. The vertical stability criterion will be $\vert a_{v}\vert<2.$ For the computation of the families of periodic orbits and their stability indexes we have performed the classical predictor-corrector scheme integrating both the equations of motion of the system and the variational equations with a integrator Adams-Bashfort-Moulton, the tolerances for the relative and absolute error are $2.22\times10^{-14}$ and $10^{-16}$ respectively.

We have already mentioned that the equations of motion (\ref{finalequations}) posses symmetries respect both axis, these symmetries have useful implications in our computations. Let $S$ be the symmetry respect to the $x-$ axis, that is, $S(x,y,\dot{x},\dot{y},t)=(x,-y,-\dot{x},\dot{y},-t)$ and let $S'$ be the symmetry respect to the $y-$ axis, that is, $S(x,y,\dot{x},\dot{y},t)=(-x,y,\dot{x},-\dot{y},-t)$, the composition of both symmetries produces a symmetry respect the origin $S\circ S'(x,y,\dot{x},\dot{y},t)=(-x,-y,-\dot{x},-\dot{y},t)$, as a consequence, we have that to any orbit there exists another orbit (possibly the same) which is symmetrical respect to the origin of coordinates, therefore, given a family of periodic orbits, there corresponds a symmetrical family. For instance, we will see that there is a family of periodic orbits around the equilibrium point $L_{1}$, because of the above mentioned symmetries, there exists a symmetric family around the equilibrium point $L_{2}$ with the same characteristics under the symmetries. Therefore, for each family explored in this work, there exists its symmetric family which is not shown in the further discussion, however, the reader must keep in mind its existence.

\subsection{The Sun-Jupiter value}

If we consider small values of the mass parameter $\mu$, we can consider the planar H4BP as a perturbation of the classical Hill's problem. More precisely, the Hamiltonian (\ref{rotatedspatialhamiltonian}) for the planar case can be written as

\begin{equation}\label{perturbedhillproblem}
H=\frac{1}{2}(p_x^2+p_y^2)+yp_x-xp_y-x^2+\frac{1}{2}y^2-\frac{1}{\sqrt{x^2+y^2}}+\frac{9}{8}\mu(x^2-y^2)+\mathcal{O}(\mu^2).
\end{equation}
 
If we neglect the higher order terms, the above Hamiltonian becomes

$$
H=H_{Hill}+\mu P(x,y),
$$

where $P(x,y)=\frac{9}{8}(x^2-y^2)$ is the main the main term of the perturbation.

In the Sun-Jupiter system we consider the value $\mu=0.00095$, so, for this particular case, we can think the system as a perturbation of the classical Hill problem therefore it is expected that the structure of the periodic orbits of the H4BP be similar to the ones of the classical Hill's problem. Our numerical explorations shows that this is the case, however, such explorations show new features regarding the stability of the families of periodic orbits.

\subsection{The family $g$}
\begin{figure}
\centering
\includegraphics[width=1.5in]{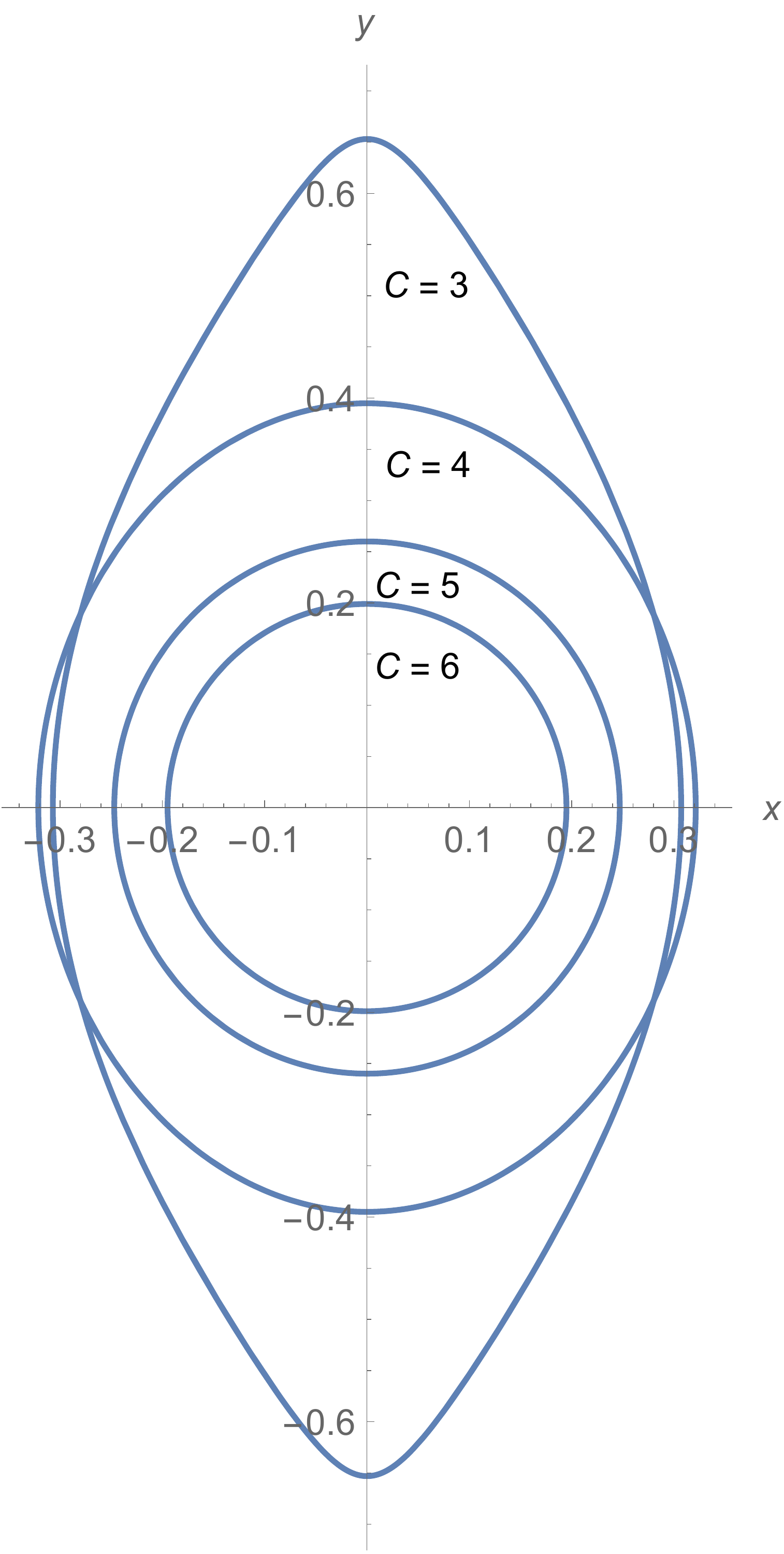}
\caption{Evolution of the family $g$.\label{g1}}
\end{figure}

\begin{figure}
\centering
\includegraphics[width=1.0in]{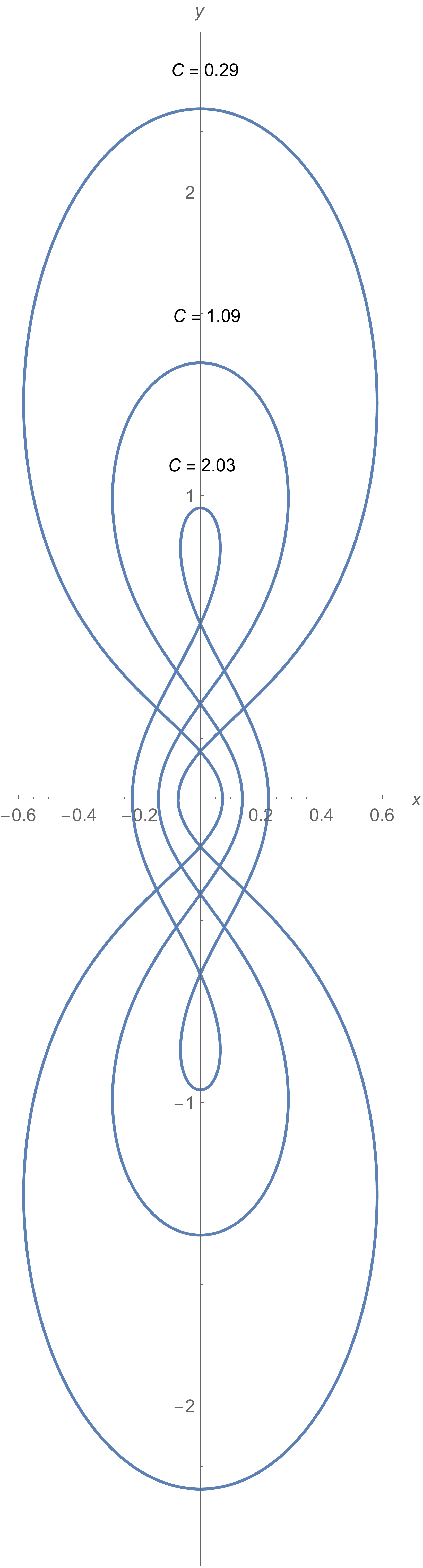}
\caption{Evolution through collision of the family $g$.\label{g2}}
\end{figure}

\begin{figure}
\centering
\includegraphics[width=1.0in]{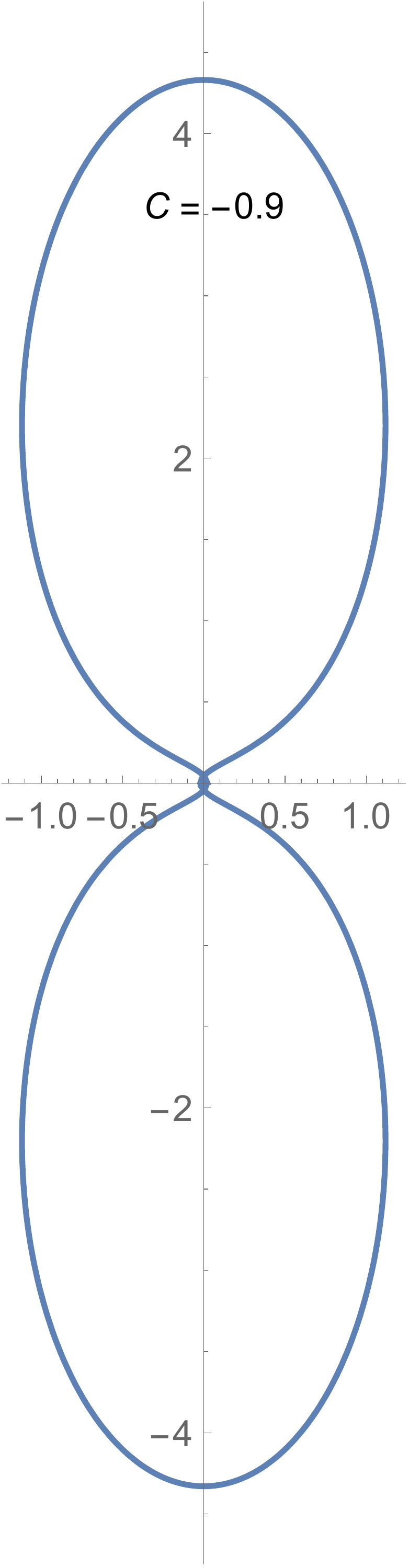}
\caption{Orbit near to collision for the family $g$.\label{gfinal}}
\end{figure}

\begin{figure}
\centering
\includegraphics[width=2.0in]{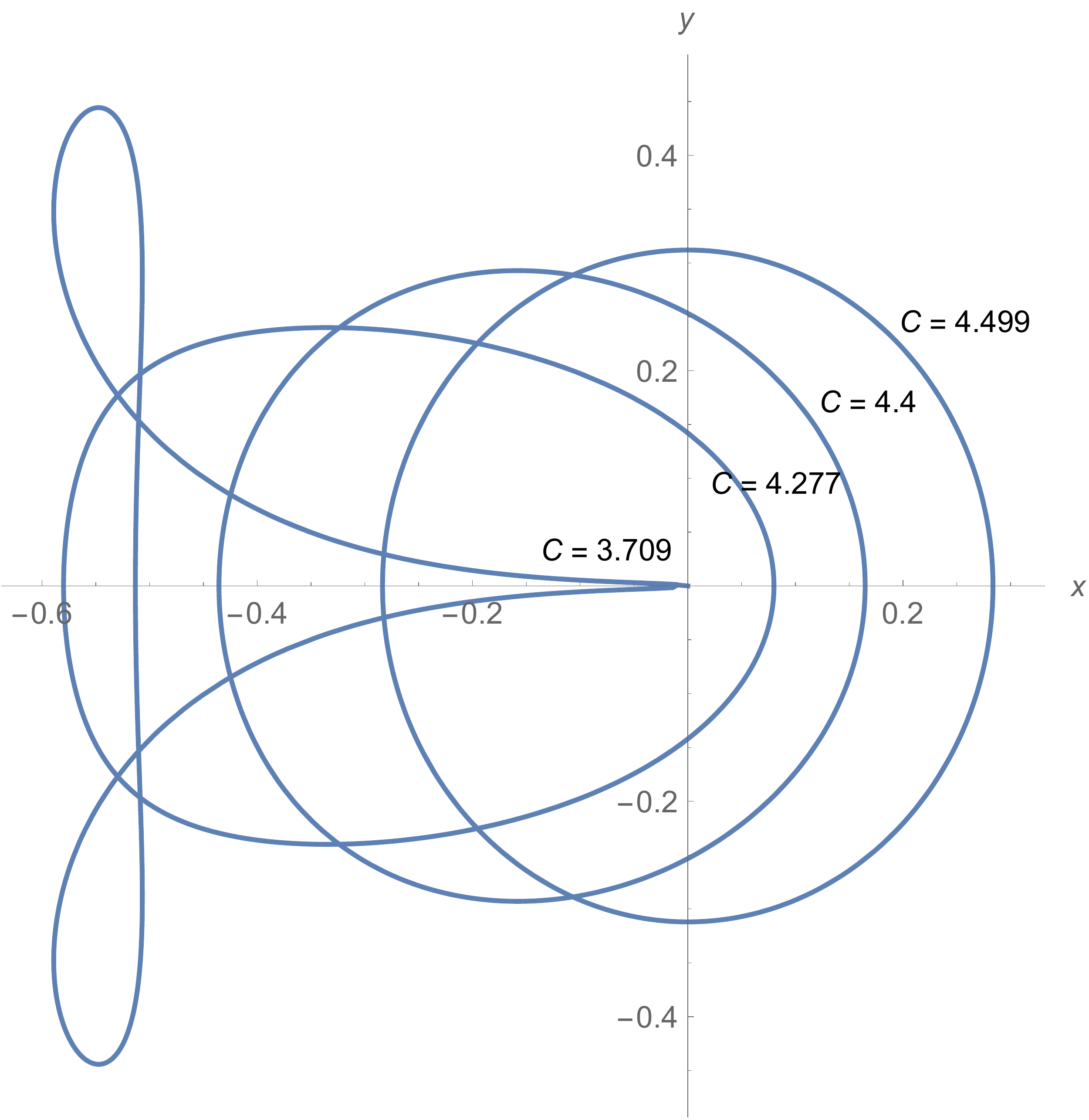}
\caption{Evolution of the lower bifurcation branch of the family g.\label{gl1}}
\end{figure}

\begin{figure}
\centering
\includegraphics[width=1.0in]{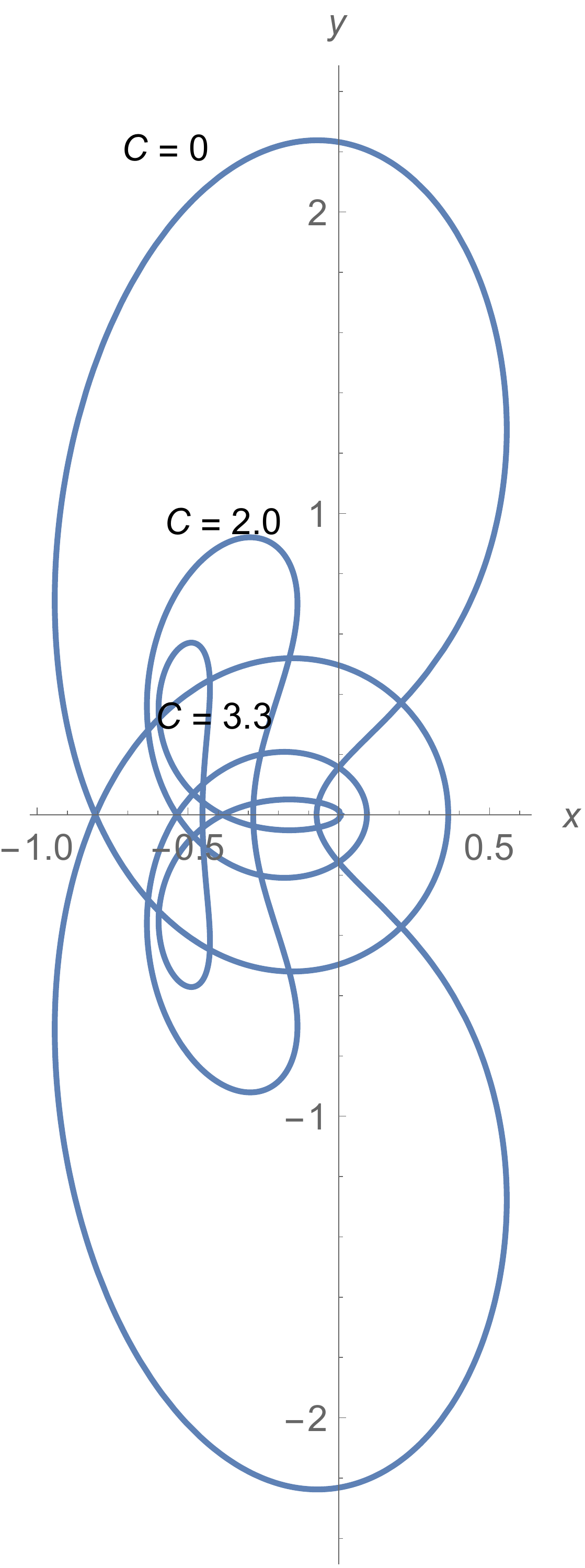}
\caption{Evolution of the lower bifurcation branch of the family g.\label{gl2}}
\end{figure}

\begin{figure}
\centering
\includegraphics[width=1.0in]{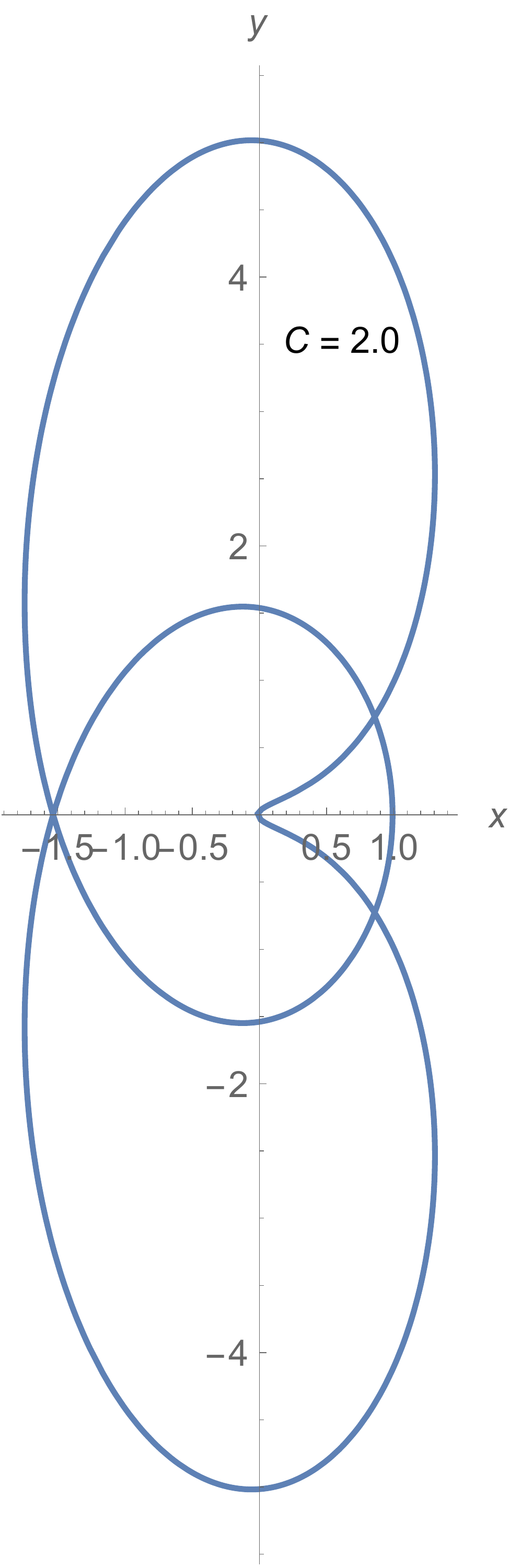}
\caption{Orbit near to collision for the lower bifurcation branch of the family g.\label{gl2final}}
\end{figure}

\begin{figure}
\centering
\includegraphics[width=2.0in]{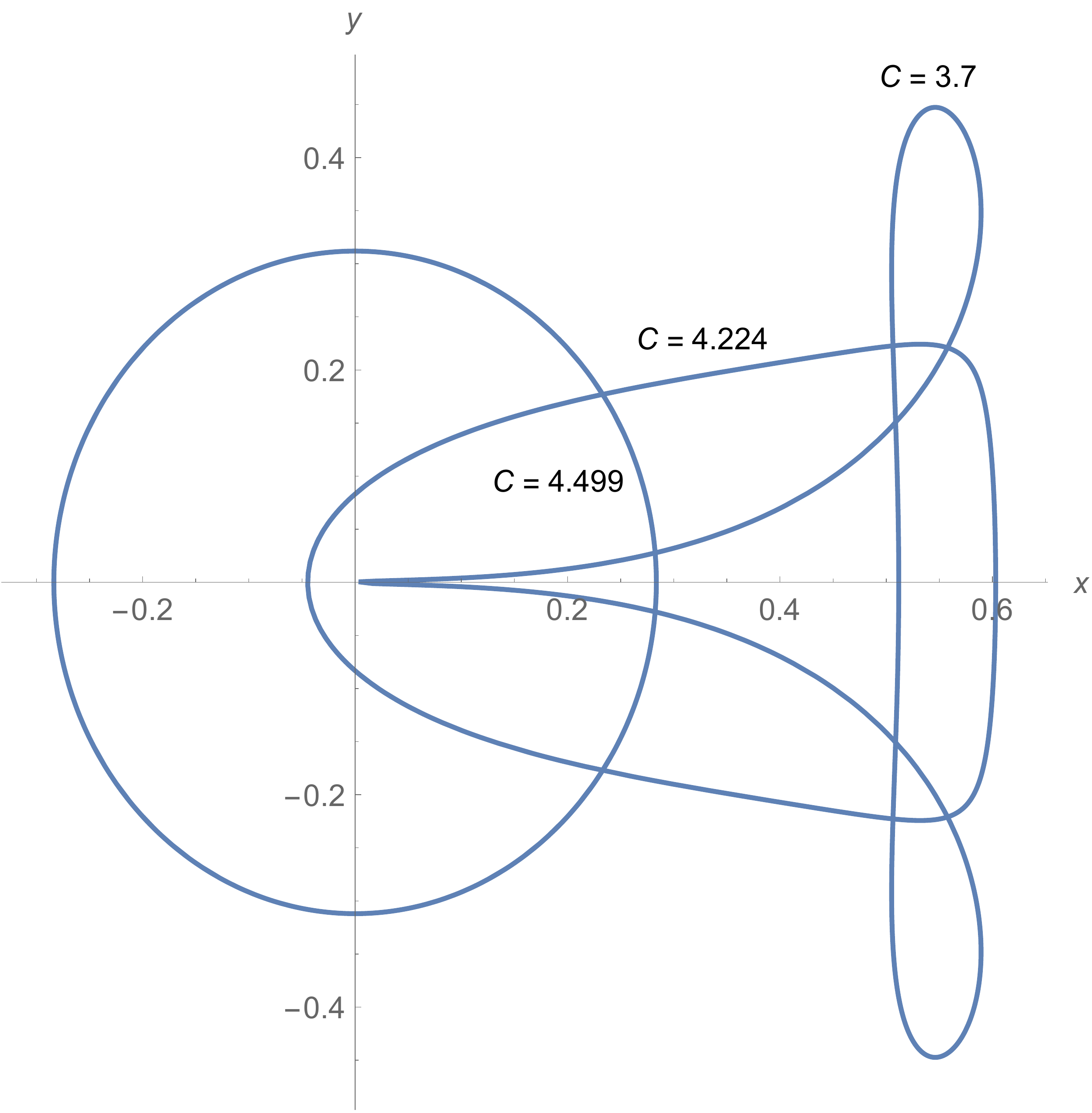}
\caption{Evolution of the upper bifurcation branch of the family g.\label{gu1}}
\end{figure}

\begin{figure}
\centering
\includegraphics[width=1.0in]{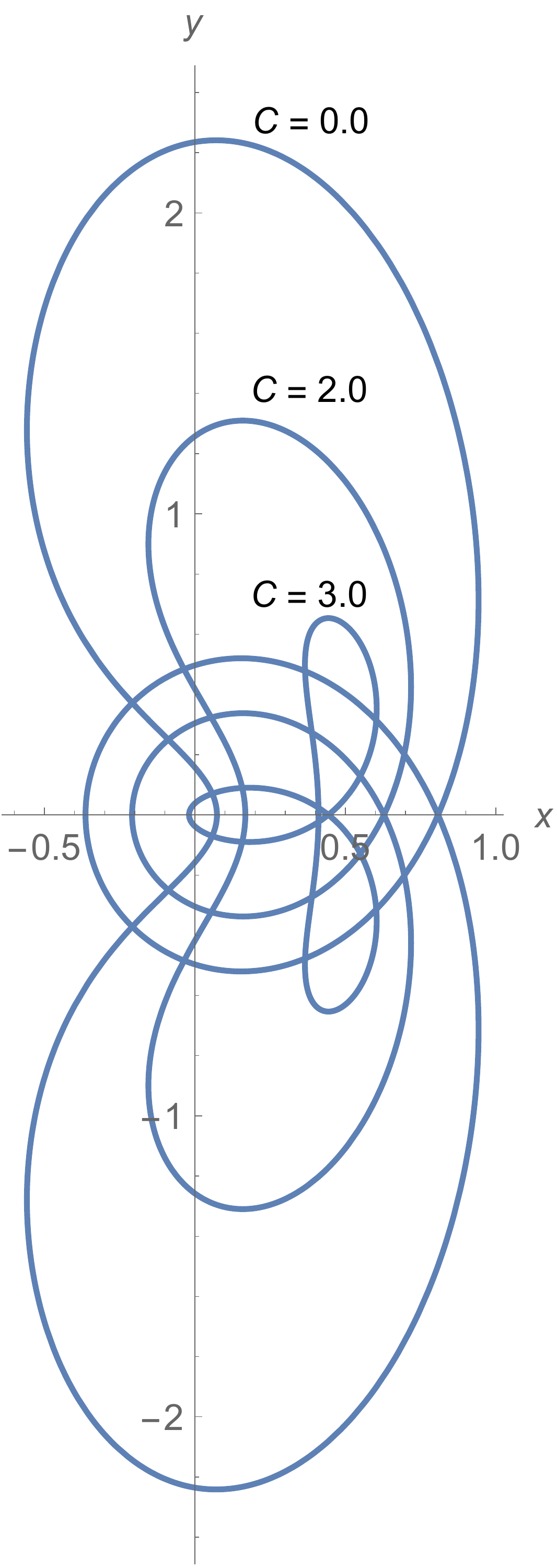}
\caption{Evolution of the upper bifurcation branch of the family g.\label{gu2}}
\end{figure}

\begin{figure}
\centering
\includegraphics[width=1.0in]{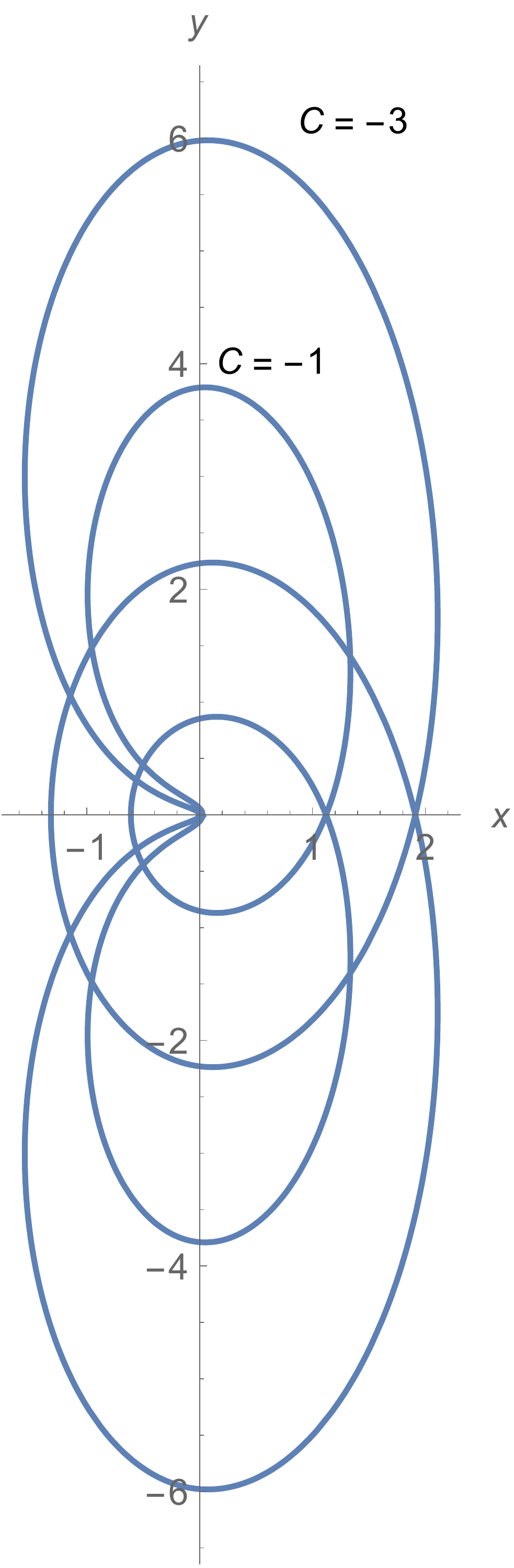}
\caption{Evolution towards collision of the lower bifurcation branch of the family g.\label{gufinal}}
\end{figure}

This family begins with infinitesimal direct circular periodic orbits around the tertiary, the evolution of the family is quite similar to the correspondent family $g$ for the classical Hill's problem, the size of the orbits increases as the Jacobi constant is decreased, the shape of the orbits are quite similar as in the Hill's problem and they tend to collision asymptotically, i.e., the initial condition $x_{0}$ tends to zero for $C\rightarrow-\infty$. See figures \ref{g1}, \ref{g2} and \ref{gfinal}. Throughout the evolution of this family, we found a bifurcation point at $C=4.4983599991.$ as a consequence we emerge two branches from this point, we have called $g$-lower and $g$-upper to the respective families. The evolution of the this branch $g$-lower is shown in the figures \ref{gl1}, \ref{gl2} and \ref{gl2final}. In the figure \ref{gl1} we can observe that as $C$ decreases from the value of bifurcation the size of the orbits increases until a collision with the tertiary is reached, after this collision, a new loop appear on the orbit, this loop and the size of the orbit increase as we keep decreasing the value of $C$, a new collision orbit is reached approximately in $C=-15.0$. The evolution of the branch $g$-upper is shown in the figures \ref{gu1}, \ref{gu2} and \ref{gufinal}. Analogously as in the lower branch, the size of the orbits increases until to reach a collision orbits as the value of $C$ decreases, throughout this evolution we can observe a change of stabilty of the orbits. After the collision orbits, a new loop appears in the orbits, this loops and the size of the orbits increase until the orbits tend asymptotically to a collision orbit similarly as in the $g$ family. The family $g$ is horizontally and vertically stable (from now on, we refer to such orbits as bistable) in the interval $(-\infty,4.4979999991]$ where $a_{h}=1$, then all of the members of this family becomes horizontally unstable as $C$ decreases monotonically, such change of stability was firstly reported in \cite{BurgosIII}. The orbits remain vertically stable until the critical value $C=3.0639999991$ where $a_{v}=-1$, then the coefficient $a_{v}$ increases monotonically until a new critical value is reached in $C=1.39$ where $a_{v}=1$. As the value of $C$ keeps decreasing, the value of $a_{v}$ increases monotonically as the family tends asymptotically to collision. In the branch $g$-upper we found a rich structure regarding the stability of the orbits. As $C$ decrases monotonically after the value of bifurcation, the orbits are bistable, however a critical value is found in $C=4.43$ where $a_{v}=-1$ then the coefficient increases and the family remains bistable until the critical value $C=4.29734$ where $a_{h}=-1$, after this value both coefficients increase and we find other critical values in  $C=4.2799999991$ where $a_{v}=1$ and $C=4.2703499991$ where $a_{h}=1$. As $C$ is decreased, the coefficient the value of $a_{h}$ increases and therefore the orbits become horizontally unstable but the vertical stability remains until the value $C=4.2$ where $a_{v}=-1$. After this value the family tends to collision and the family becomes highly unstable. The members of the branch $g$-lower are bistable after the bifurcation point and both the coefficients decrease their value monotonically until a critical value is reached in $C=4.4359999991$ where $a_{v}=-1$, after this value the value of $a_{h}$ keeps decreasing monotonically but the coefficient $a_{v}$ begins to increase its value until a new critical value is reached in $C=4.2841$ where $a_{v}=1$, now the family tend to collision and the coefficient $a_{h}$ increases its value and a critical orbit is reached in $C=4.27034$ where $a_{h}=1$, after this value such coefficient increases as $C$ decreases monotonically and therefore the family becomes unstable.

\subsection{The family $f$}

\begin{figure}
\centering
\includegraphics[width=1.5in]{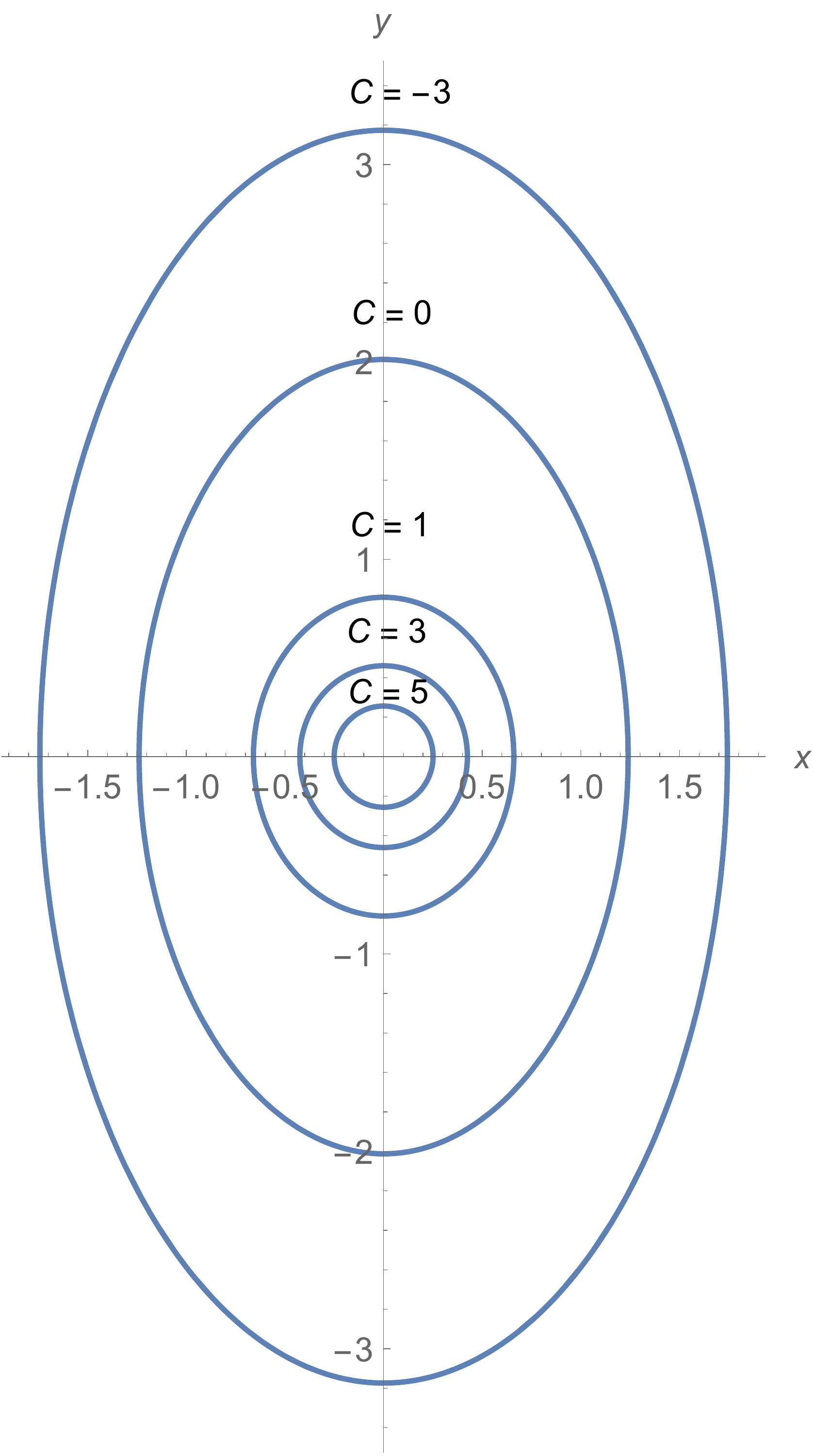}
\caption{The evolution of the family $f$ of retrograde periodic orbits.\label{f}}
\end{figure}

The members of this family are retrograde orbits around the tertiary. The evolution of this family is trivial, the size of the orbits increases as the value of $C$ decreases, we observe that this behaviour remains with no changes as $C$ decreases indefinitely. All of the members of this family are bistable. See figure \ref{f}.

\subsection{The family $a$}

\begin{figure}
\centering
\includegraphics[width=1.0in]{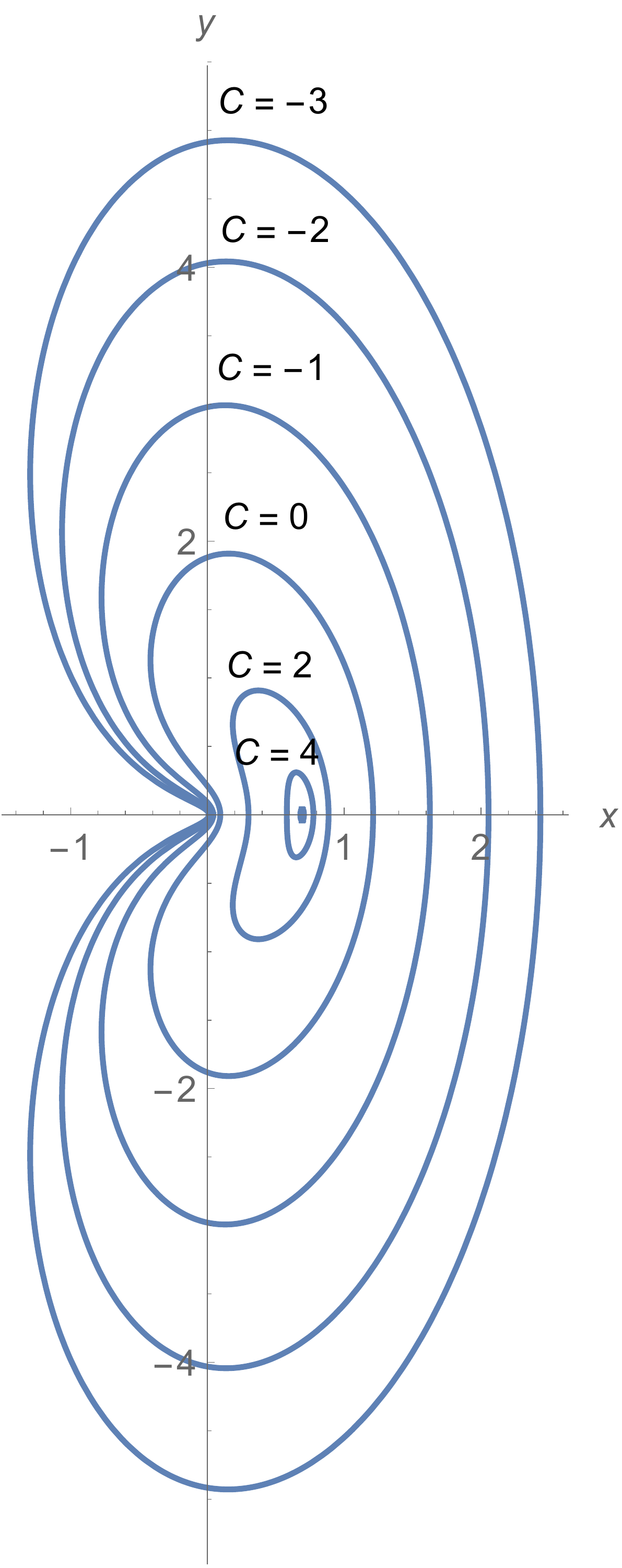}
\caption{The evolution of the family of horizontal Lyapunov orbits around the equilibrium point $L_{1}$.\label{a}}
\end{figure}

\begin{figure}
\centering
\includegraphics[width=2.0in]{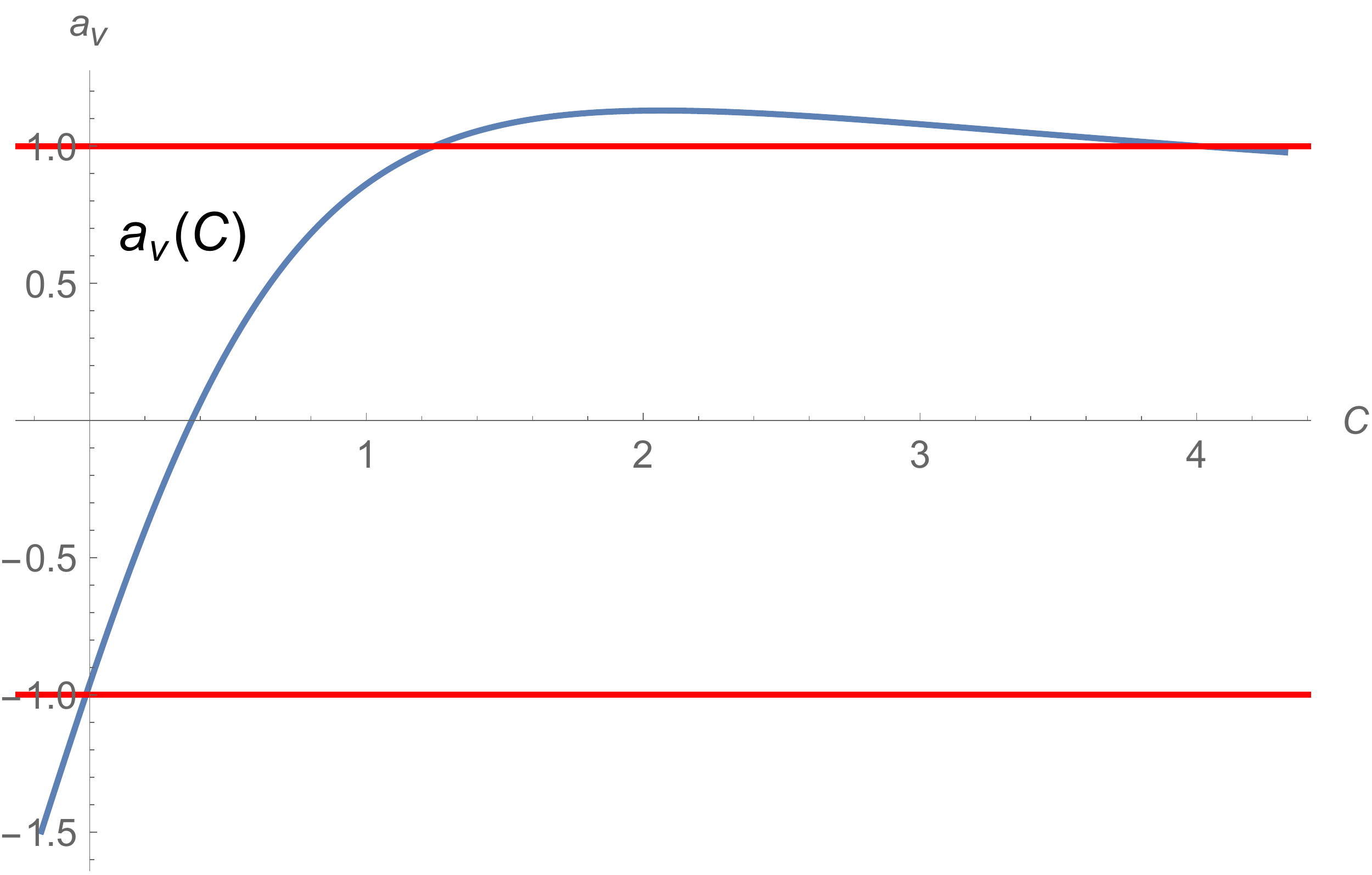}
\caption{Evolution of the coefficient of vertical stability for the family $a$.\label{verticala}}
\end{figure}

The orbits of this family come from the horizontal Lyapunov orbits around the equilibrium point $L_{1}$. As the value of $C$ decreases without bound, the size of the orbits increases and the family tends asymptotically to collision. The orbits are highly horizontal unstable, however, there are several critical orbits where the coefficient of vertical stability changes from stable to unstable type. The values of $C$ for the critical vertical orbits are $C=4.0059999991$, $C=1.246$ and $C=-0.0129999991$, see figure \ref{verticala}. It is worth mentioning that for the classical Hill's problem, a similar behaviour was reported \cite{Gomez} for the vertical stability of the Liapunov orbits. The behaviour of the symmetric family of periodic of horizontal Lyapunov orbits around the equilibrium point $L_{2}$ is completely analogous due to the symmetries of the equations. See figure \ref{a}.

\subsection{The family $H_{b}$}

\begin{figure}
\centering
\includegraphics[width=3.5in]{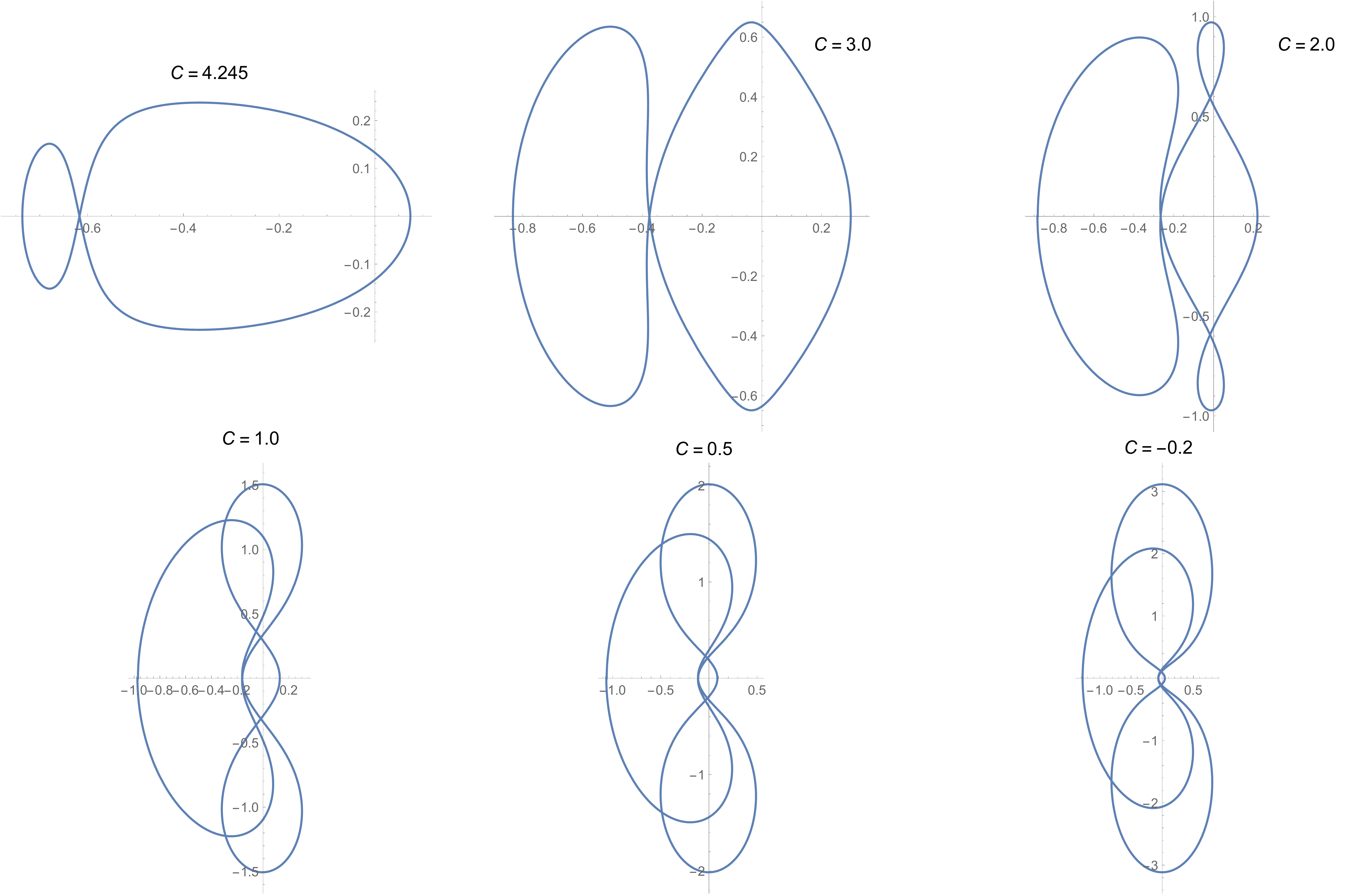}
\caption{The evolution of the family $H_{b}$.\label{familyg2}}
\end{figure}

\begin{figure}
\centering
\includegraphics[width=2.5in]{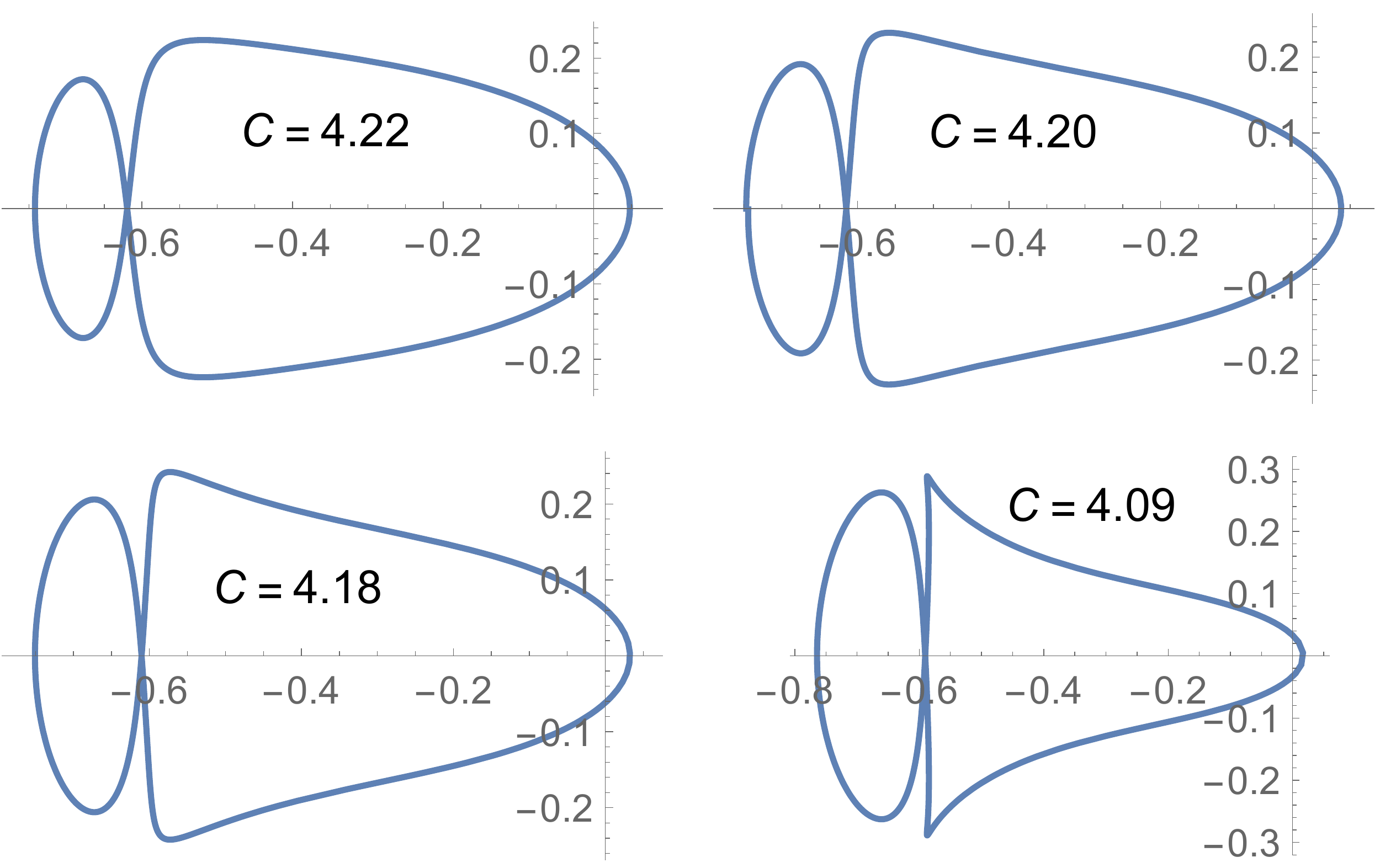}
\caption{The evolution of the second branch of the family $H_{b}$.\label{g2b2}}
\end{figure}

\begin{figure}
\centering
\includegraphics[width=2.0in]{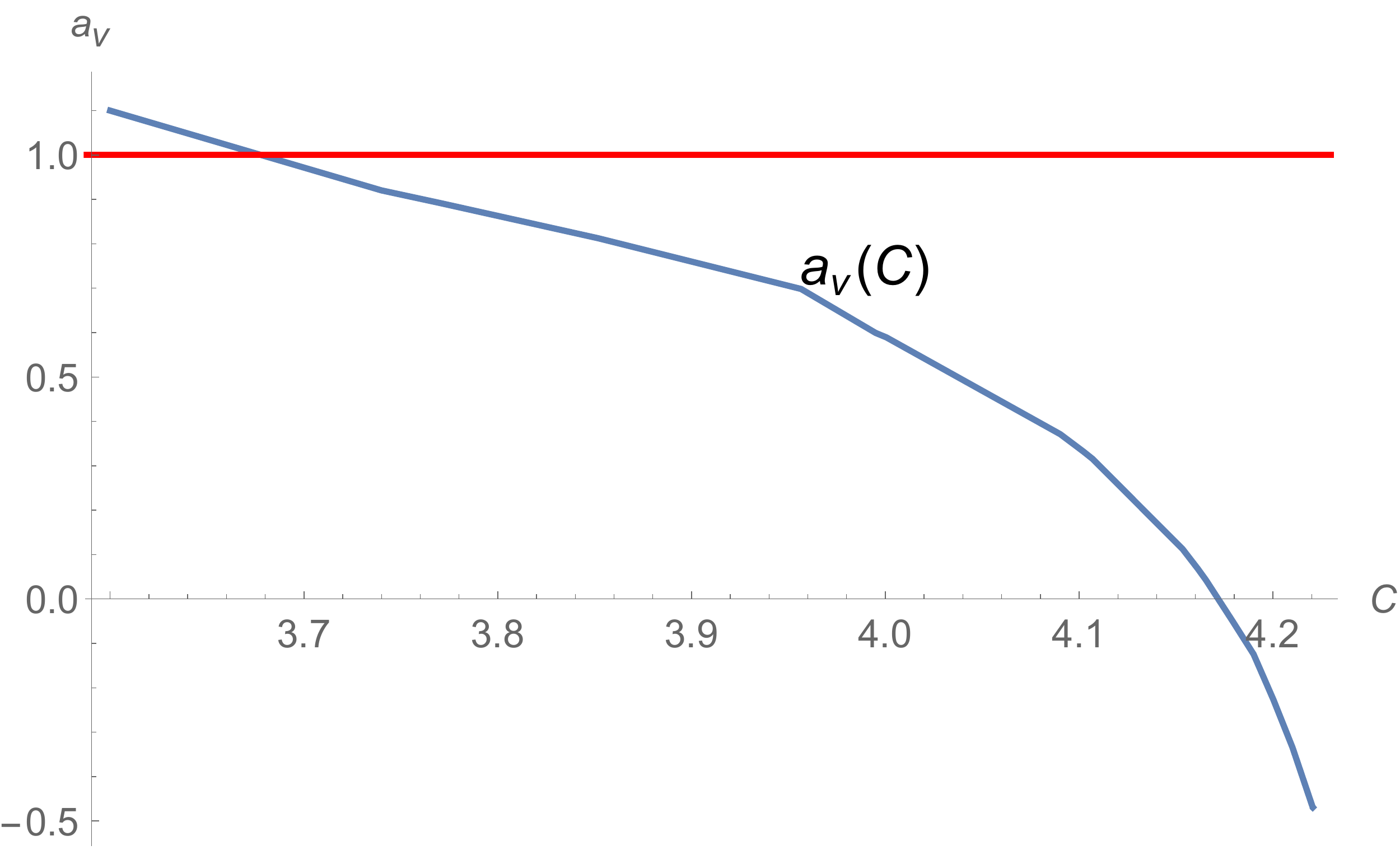}
\caption{Evolution of the coefficient of vertical stability for the second branch of the family $H_{b}$.\label{stabilityg2b2}}
\end{figure}

\begin{figure}
\centering
\includegraphics[width=2.0in]{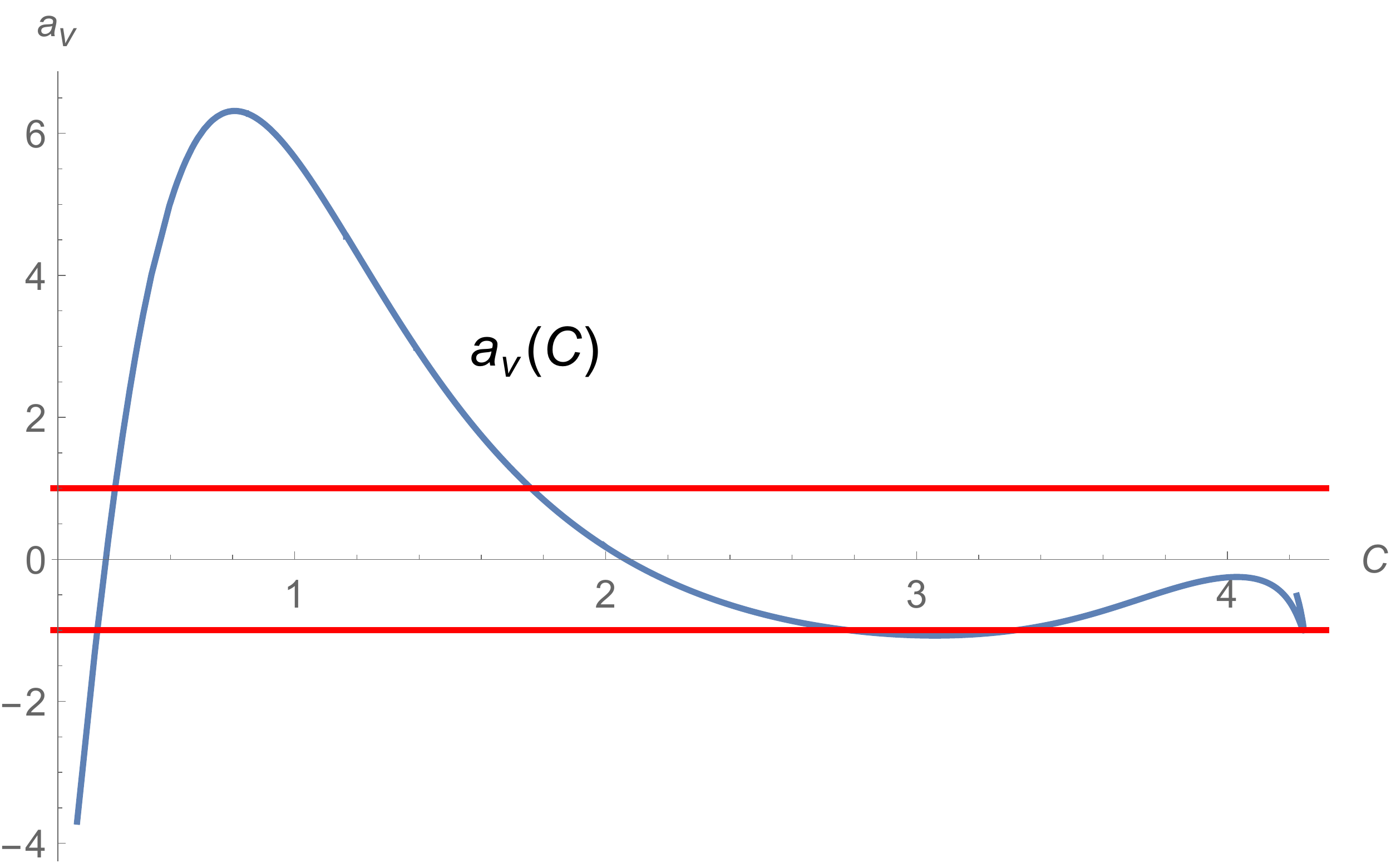}
\caption{Evolution of the coefficient of vertical stability for the family $H_{b}$.\label{verticalg2}}
\end{figure}

We have named $H_{b}$ to this family because is analogous to the that family  studied in \cite{Henon2003}. The family has a maximum in $C=4.2450522241$ and two branches emerge form this point as the value of $C$ decreases. We show the behaviour of one of the branches in the figure \ref{familyg2}. The orbits increase its size as $C$ decreases monotonically, we observe that the family tends to collision as we approach to the value $C=-0.2016399991$.
When the orbits belonging to the first branch are near to collision, they are highly vertical and horizontal unstable, in fact all of the orbits of this branch are horizontal unstable, however, there are several critical orbits where the coefficient of vertical stability changes from stable to unstable type, see figure \ref{verticalg2}. The values of $C$ for the critical vertical orbits are $C=0.3649$ where $a_{v}=-1$, $C=0.422$ where $a_{v}=1$, $C=1.759999991$ where $a_{v}=1$, $C=2.7859999991$ where $a_{v}=-1$, $C=3.3159999991$ where $a_{v}=-1$. The behaviour of the orbits belonging to the second branch is shown in the figure \ref{g2b2},  all of the orbits are horizontally unstable, however, the orbits change from vertically stable to unstable, see figure \ref{stabilityg2b2}.

\subsection{The family $H_{a}$} 
\begin{figure}
\centering
\includegraphics[width=2.0in]{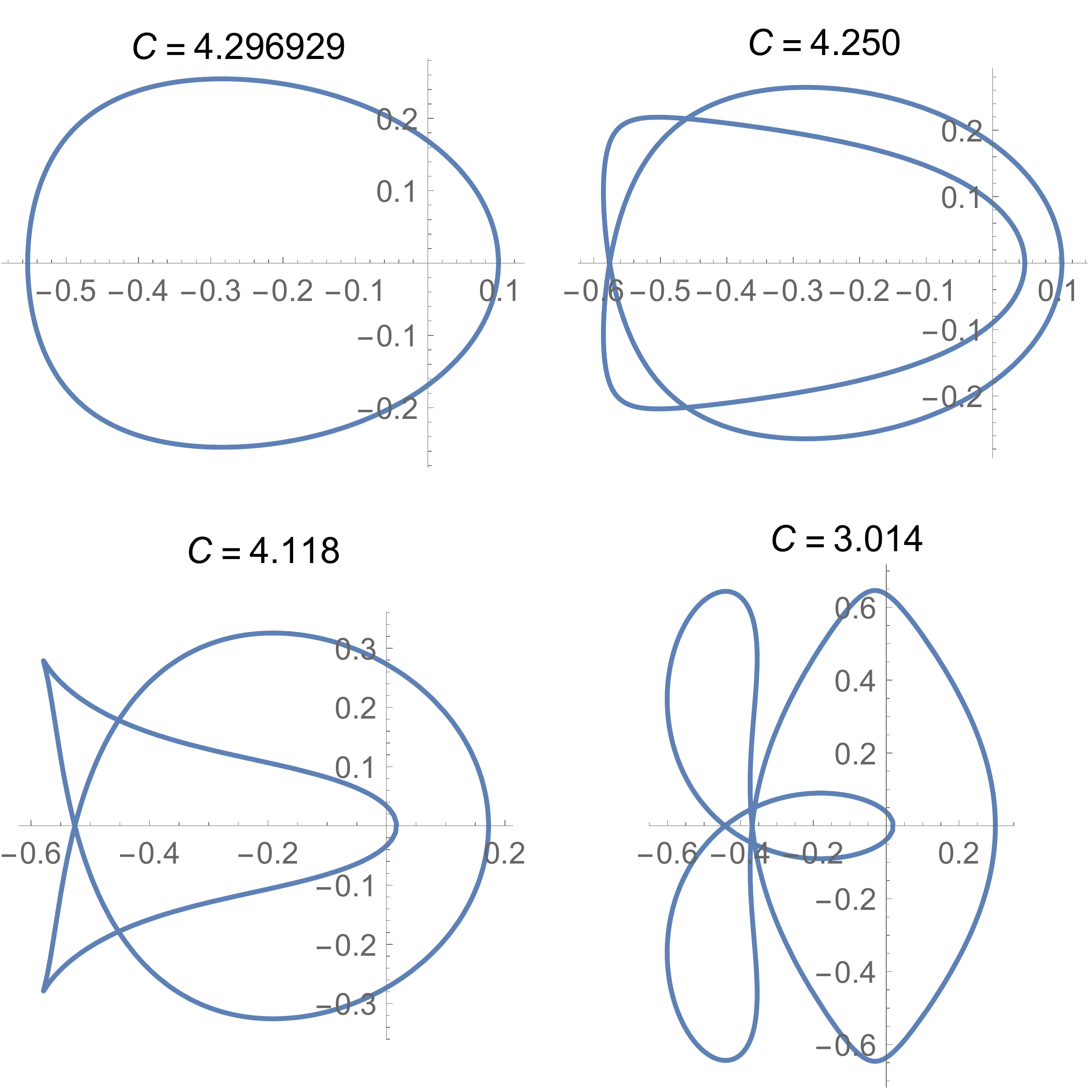}
\caption{The evolution of the family $Ha$.\label{ha}}
\end{figure}

This family emerges as a bifurcation from the family $g-$upper in the bifurcation point $C=4.1178$ where $a_{h}=-1$. As $C$ increases monotonically, we observe that the family tends to collision, see figure \ref{ha}. If we decrease the value of $C$ monotonically, we note that both coefficients $a_{h}$ and $a_{v}$ grow monotonically, in $C=4.2841299991$ we have a critical vertical orbit with $a_{v}=1$, after this point, $a_{h}$ continues increasing but $a_{v}$ starts decreasing and in $C=4.2795499991$ we have another critical orbit with $a_{v}=1$, in the range $[4.2795499991,4.2672499991)$ the family remains bistable and a new critical horizontal orbit is reached in $C=4.2672499991$ where $a_{h}=-1$, after this value, $a_{h}$ decreases monotonically and the family becomes horizontally unstable, however, the family remains vertically stable until the value $C=4.1208$ where $a_{v}=1$ then the family becomes vertically unstable.

\subsection{The short and long period families}
\begin{figure}
\centering
\includegraphics[width=1.2in]{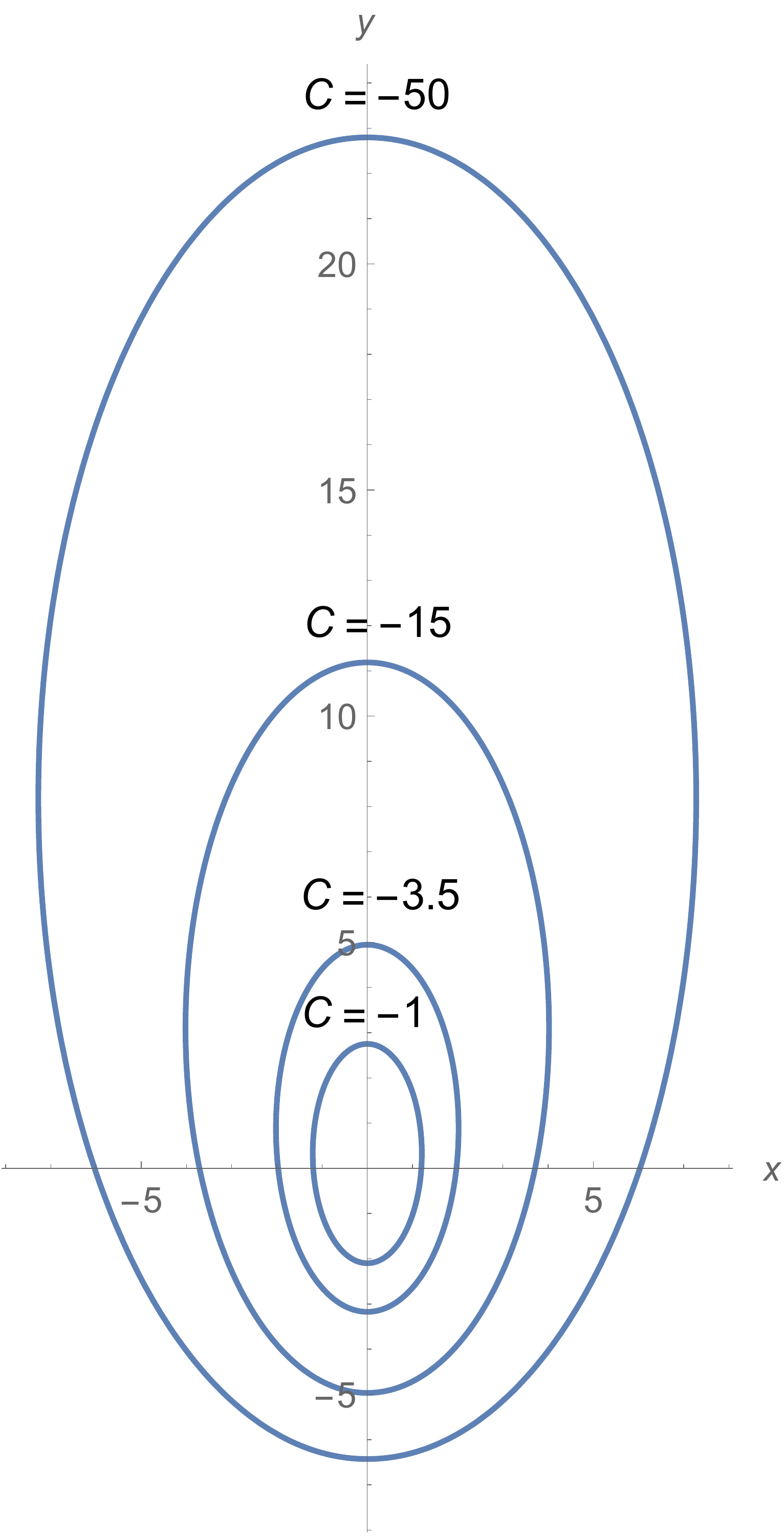}
\caption{Evolution of the short period family around $L_{3}$.\label{short}}
\end{figure}

As a consequence of the presence of a second massive body, we have two equilibrium points whose linear stability was studied in the section \ref{lineardynamics}. Such analysis suggest the existence of the two families of Lyapunov orbits, the so called short and long period families. In the classical Hill's problem there are not analogous families of such Lyapunov orbits because the equilibrium points $L_{3}$ and $L_{4}$ do not exist in the classical problem, however, it is well known that in the R3BP families of such orbits exists around the so called equilibrium points $L_{4}$ and $L_{5}$, if the value of the mass parameter is less than the Routh's critical value \cite{Sz}. In the classical work \cite{Deprit3}, we can find a deep analysis of the behaviour of the orbits obtained by numerical continuation of the infinitesimal Lyapunov orbits around $L_{5}$ and the evolution of such families is shown as the Jacobi constant changes its value. It is worth mentioning that the evolution of analogous families in our case inherits some behavior of the R3BP and the full R4BP \cite{PapaIII}. The short period family is obtained by the numerical continuation of the short period ellipses produced by the linear dynamics studied in the section \ref{lineardynamics} and the origin of coordinates has been translated to the equilibrium point. As we decrease the value of $C$ monotonically from the value $C_{L_{3}}$ the size of the ellipses increases and the shape becomes non-symmetric respect the $x$-axis, see figure \ref{short}. We observe that as $C\rightarrow-\infty$ the orbits increases its size with no bound, it is worth mentioning that the the period of the orbits also increases but very slowly. In the table \ref{valuesshortfamily} we find some initial conditions for the short period family, we have omitted the inital values of $y$ and $\dot{y}$ because the first one is zero and $\dot{y}$ is obtained from the Jacobi first integral.

\begin{table}[ht]
\caption{Some initial conditions for the short period family} 
\centering 
\begin{tabular}{c c c c}
\hline\hline 
$C$ & $x_{0}$ & $\dot{x}_{0}$ & $T$ \\ [0.5ex] 
\hline 
0.386390 & 0.0052630577 & -0.000003114  & 6.352714861 \\ 
0.000490 & 0.6349317173  & -0.0452963681  & 6.352729416 \\
-6.033910 & 2.4990322956  & -0.6970916287 & 6.352966358 \\
-99.90891 & 7.7351384429 & -6.6652904069  & 6.3561117178 \\ [1ex] 
\hline 
\end{tabular}
\label{valuesshortfamily} 
\end{table}

\begin{figure}
\centering
\includegraphics[width=0.5in]{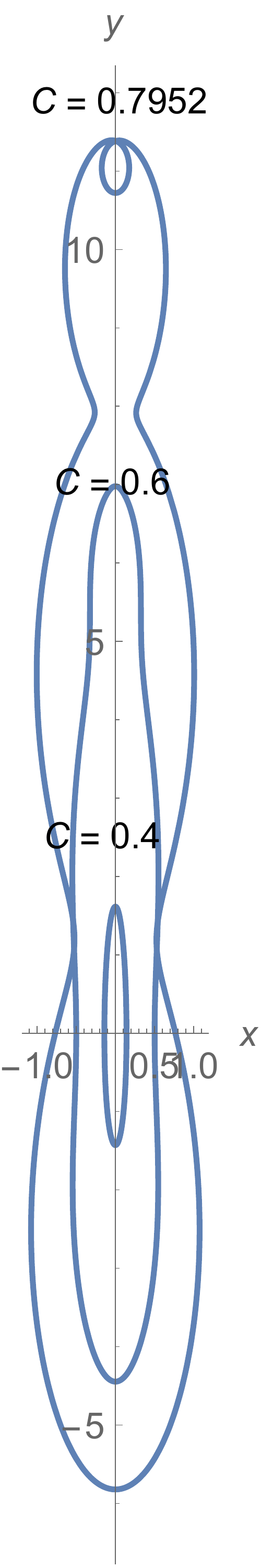}
\caption{Evolution of the first branch of the long period family around $L_{3}$.\label{longphase1}}
\end{figure}

\begin{figure}
\centering
\includegraphics[width=2.0in]{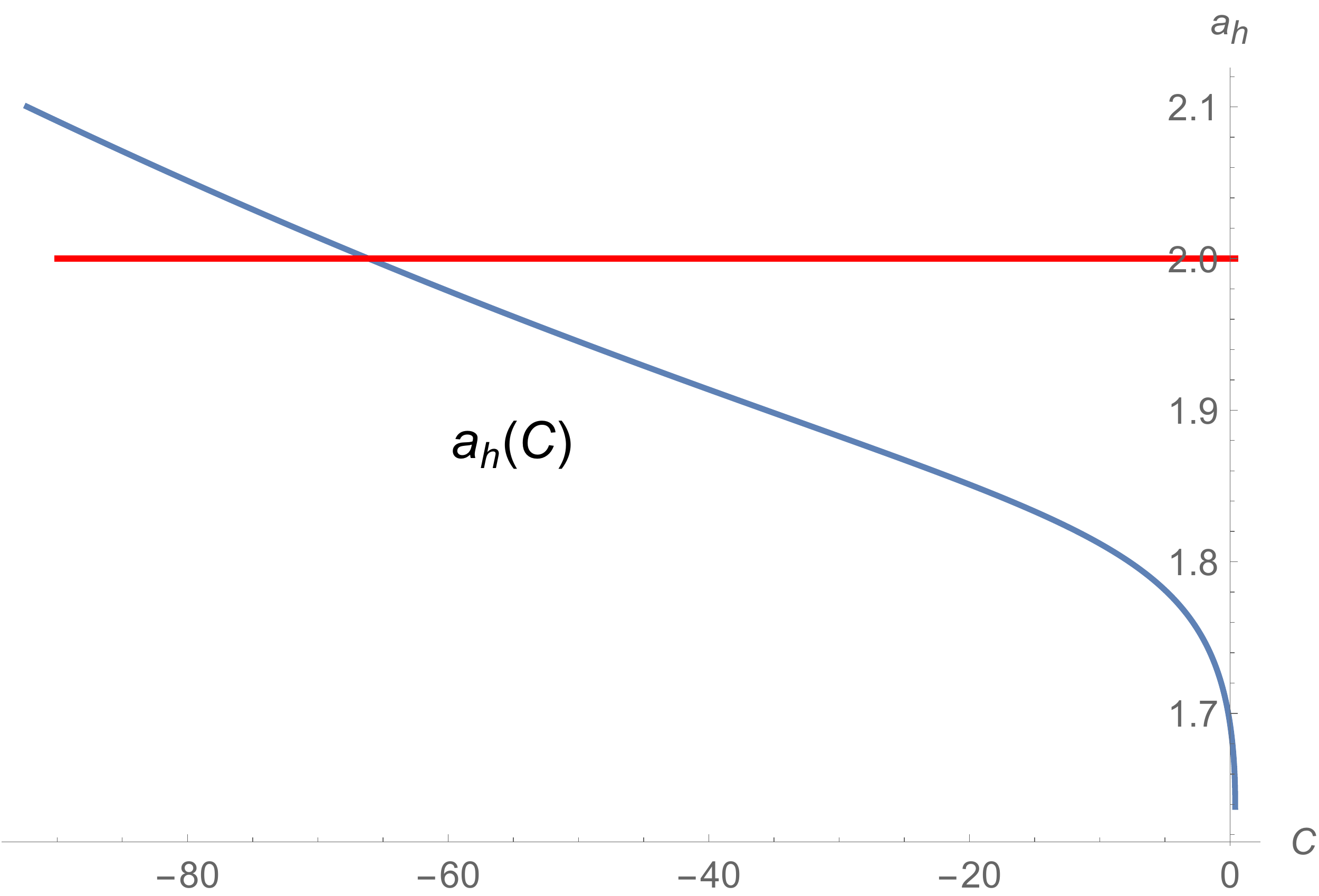}
\caption{Evolution of the horizontal stability curve of the short period family.\label{hstabilityshort}}
\end{figure}

\begin{figure}
\centering
\includegraphics[width=2.0in]{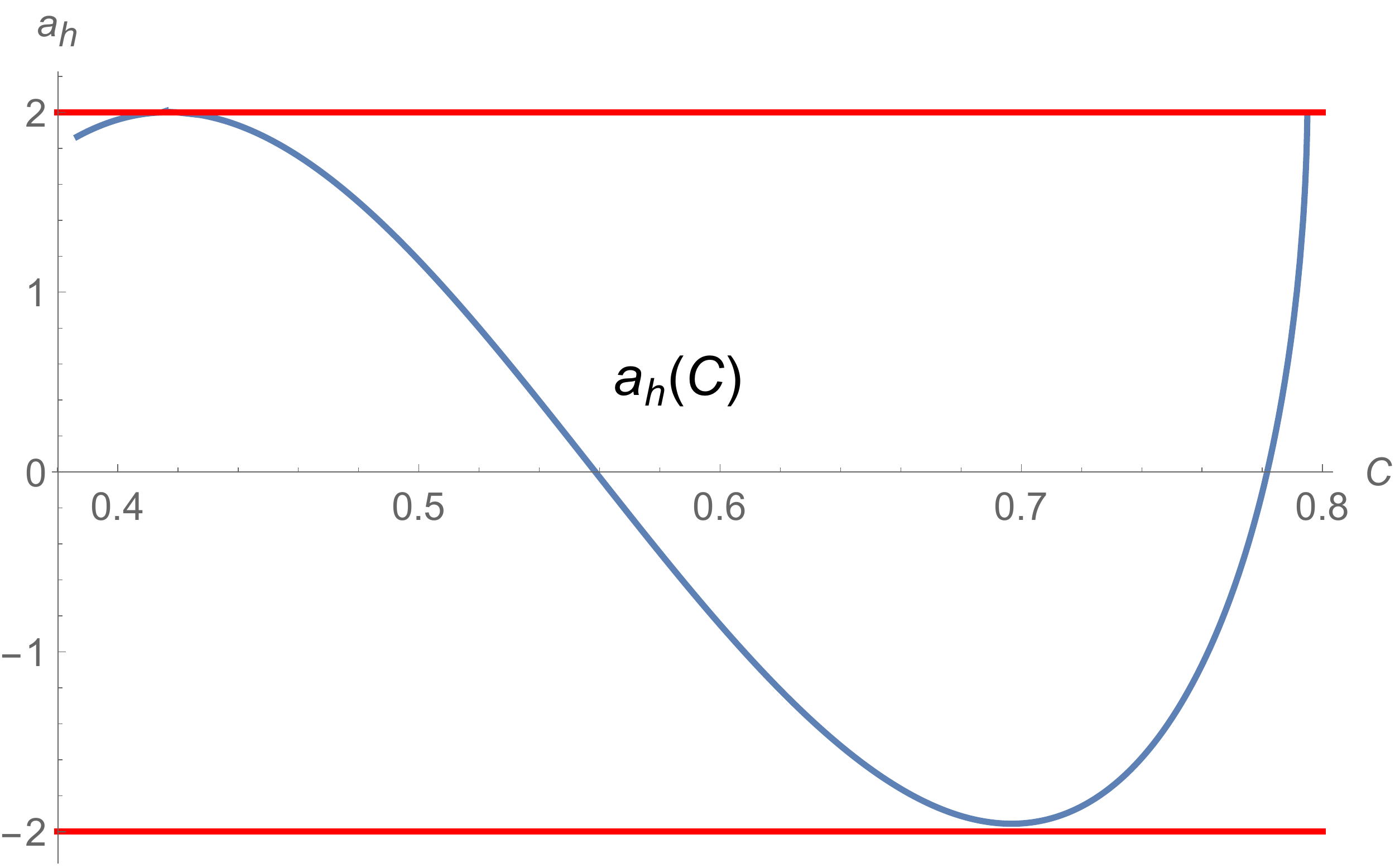}
\caption{Evolution of the horizontal stability curve of the long period family.\label{hstabilitylong}}
\end{figure}

\begin{figure}
\centering
\includegraphics[width=2.0in]{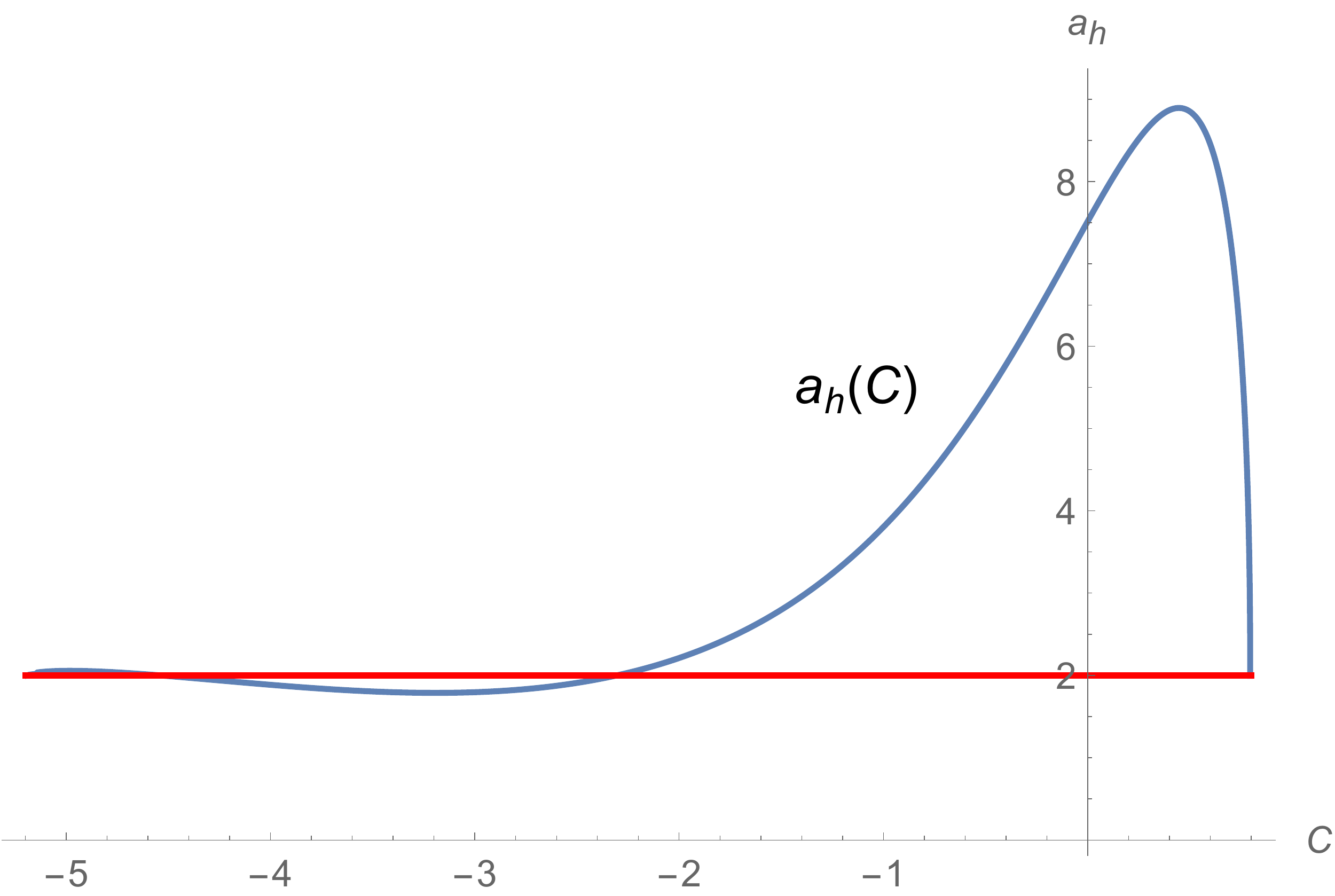}
\caption{Evolution of the horizontal stability curve after the turning point of the long period family.\label{hstabilitylongafterfold}}
\end{figure}

\begin{figure}
\centering
\includegraphics[width=2.0in]{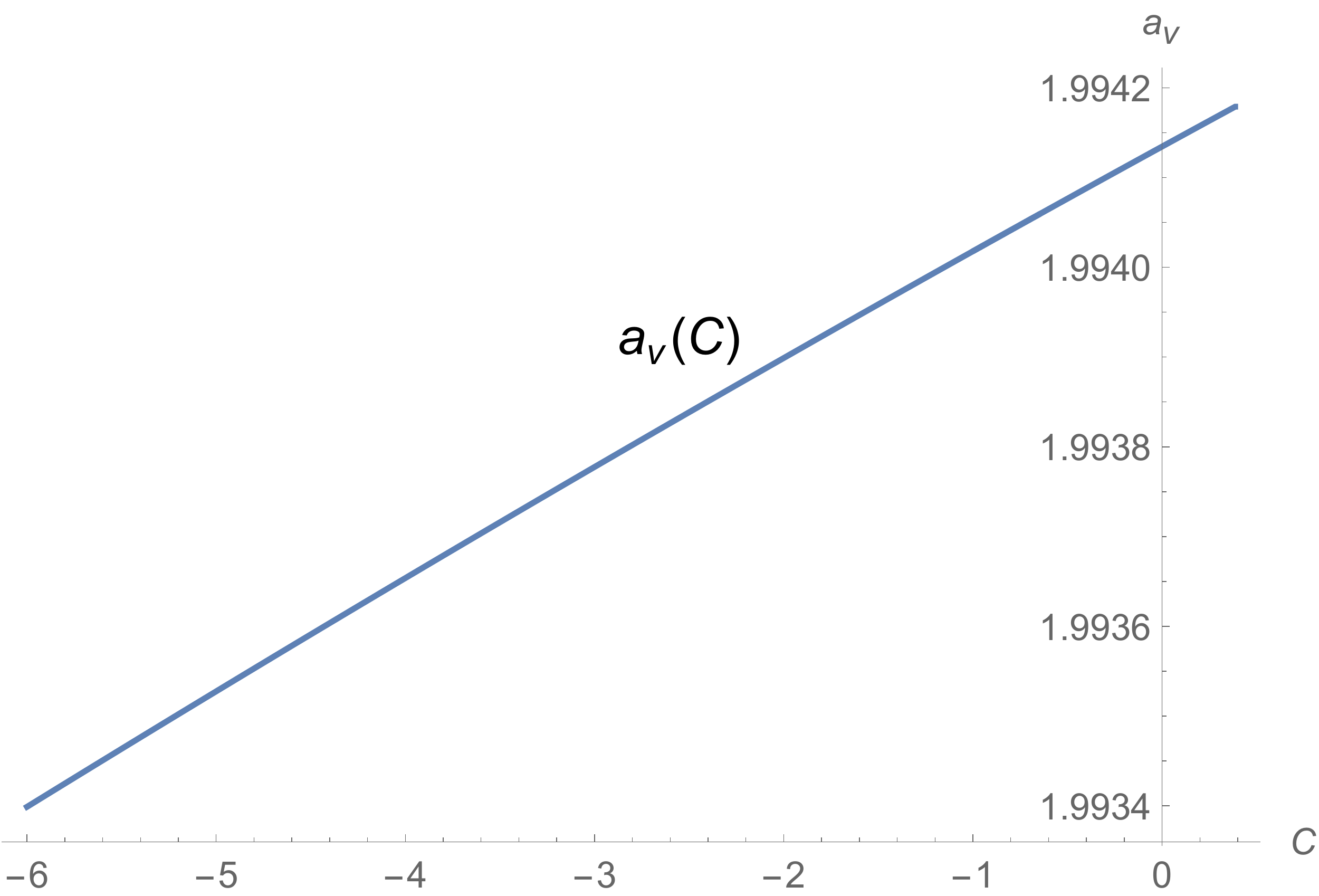}
\caption{Evolution of the vertical stability curve of the short period family.\label{verticalshort}}
\end{figure}

\begin{figure}
\centering
\includegraphics[width=2.0in]{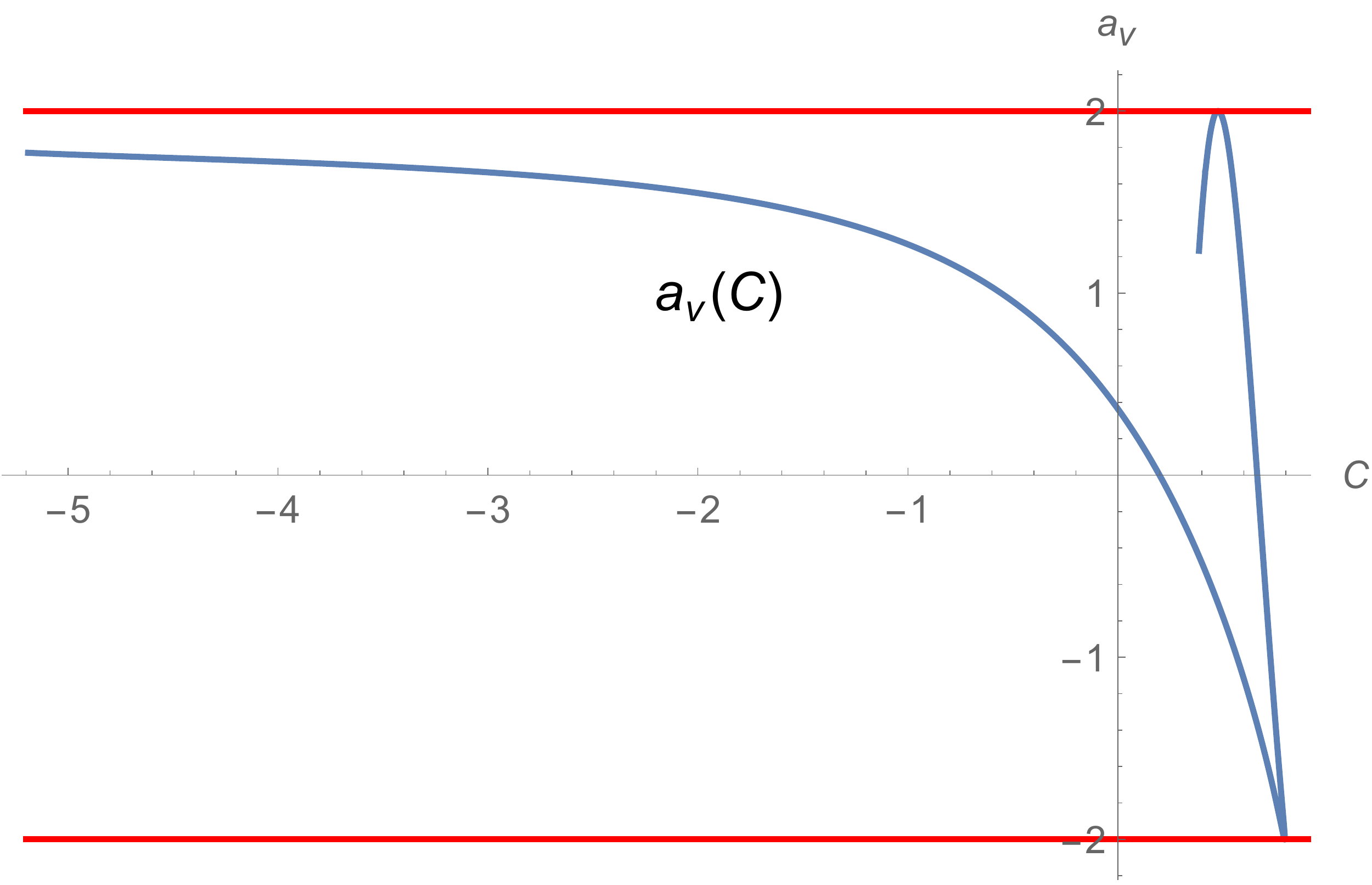}
\caption{Evolution of the vertical stability curve of the long period family.\label{verticalong}}
\end{figure}

The structure of the long period family is much more interesting, as we increase the value of $C$ the size of the orbits increase as it is shown in the figure \ref{longphase1} , in $C=0.7952$ we have a critical orbits of first type, more precisely a turning point. After this point the value of $C$ increases and the size of the orbits starts to decrease as it is shown in the figure \ref{longphase1}. It is easy to see that this branch tends to a new bifurcation point at $C=-5.1901399991$ where we obtain a orbit such that the 6 loops coincide. In the figures \ref{hstabilityshort},\ref{hstabilitylong} and \ref{hstabilitylongafterfold} we can observe the behavior of the horizontal stability of the families of periodic orbits, all of the members of the short period family are horizontally stable until the value $C=-66.109983619$ where $a_{h}:=a+d=2$, after that value all of the orbits become horizontally unstable. The stability of the second branch of the long period shows a remarkable behaviour, in the figure \ref{hstabilitylongafterfold}. We observe that the horizontal stability curve passes through the stable to the unstable zone several times, the values of $C$ for which this occurs are $C=0.415800029$, $C=0.419690042$, $C=0.7949799991$, the last value corresponds to the turning point of the branch. After this last point, all of the orbits become horizontally unstable until the value $C=-2.30591999991$, then we observe crossings between the stable and unstable areas, the values where this occurs are $C=-4.540139999$ and $C=-5.190139999$, this last value corresponds to the bifurcation with the short family. The vertical stability of both families is shown in the figures \ref{verticalshort} and \ref{verticalong}, we can observe that all the short family is vertically stable although the coefficient $a_{v}$ is very near to the critical value $a_{v}=2$, however, we can find a vertical critical value $a_{v}=-2$ in the long period family, the value of the Jacobi constant for such orbit is $C=0.794679999$ that is near to the value of $C$ for which a turning point occurs, besides this point, all of the family is vertically stable.

\begin{figure}
\centering
\includegraphics[width=1.8in]{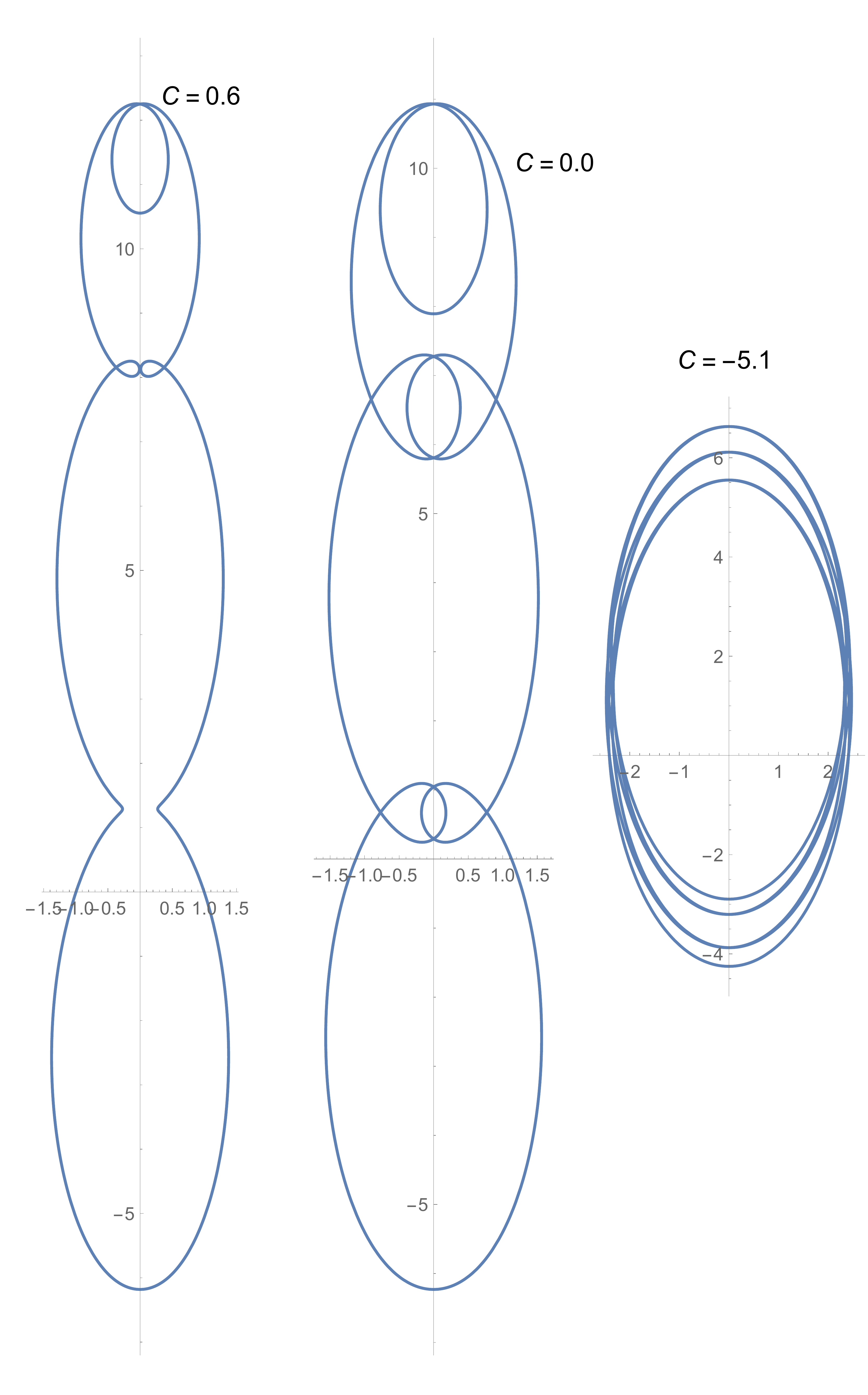}
\caption{Evolution of the first branch of the long period family around $L_{3}$.\label{longphase2}}
\end{figure}

\end{document}